\documentclass[12pt,fleqn, a4]{article}

\usepackage[french]{babel}
\usepackage{inputenc}
\usepackage[T1]{fontenc}

\usepackage{amssymb}
\usepackage{latexsym}
\usepackage{graphics}

\usepackage{graphics}
\usepackage{color}
\newcommand{\colval}{0.3}
\definecolor{colone}{gray}{\colval}

\newcommand{\dcb}{\begin{array}{lll}}
\newcommand{\dce}{\end{array}}
\newcommand{\ebe}{\begin{enumerate}\setlength{\baselineskip}{13pt}\setlength{\parskip}{5pt}}
\newcommand{\dbe}{\end{enumerate}}

\newcommand{\ibegin}{\begin{itemize}\setlength{\baselineskip}{19pt}\setlength{\parskip}{7pt}}
\newcommand{\iend}{\end{itemize}}
\newcommand{\ok}{\rule{4pt}{6pt}}
\newcommand{\desb}{\begin{description}}
\newcommand{\dese}{\end{description}}

\newtheorem{Thm}{Theorem}[section]
\newtheorem {Cor}[Thm]{Corollary}
\newtheorem {definition}[Thm]{Definition}
\newtheorem {pro}{Proposition}[Thm]
\newtheorem {Lemma}[Thm]{Lemma}
\newtheorem {rem}[Thm]{Remark}

\newtheorem {assumption}[Thm]{Assumption}

\newcommand {\bd}{\begin{definition}}
\newcommand {\ed}{\end{definition}}
\newcommand {\bpro}{\begin{pro}}
\newcommand {\epro}{\end{pro}}
\newcommand {\bl}{\begin{Lemma}}
\newcommand {\el}{\end{Lemma}}
\newcommand {\bcor}{\begin{Cor}}
\newcommand {\ecor}{\end{Cor}}
\newcommand {\brem }{\begin{rem} \rm }
\newcommand {\erem }{\end{rem}}
\newcommand{\bethe}{\begin{Thm}}
\newcommand{\ethe}{\end{Thm}}
\newcommand {\bassumption}{\begin{assumption}}
\newcommand {\eassumption}{\end{assumption}}

\def \ind{1\!\!1}
\def\cro#1{\langle #1\rangle}

\begin{document}

\begin{center}
\textbf{\Large Random time with differentiable conditional distribution function}\footnote{A preliminary version}
\end{center}

\begin{center}
Shiqi Song

{\footnotesize Laboratoire Analyse et Probabilités\\
Université d'Evry Val D'Essonne, France\\
shiqi.song@univ-evry.fr}
\end{center}

\section{Introduction}

In \cite{songmodel} a particular class of one-default market models was presented, where the default times were defined by stochastic differential equations. We learned from this study that random times $\tau$ in this class might have their conditional distribution function $u\rightarrow\mathbb{Q}[\tau\leq u|\mathcal{F}_t]$ differentiable with respect to an $\mathbb{F}$ adapted non decreasing process $A$, and the derivatives were computed in term of the stochastic flow associated with the stochastic differential equations. Various consequences of these random times were discussed based on that differentiability with respect to $A$. It then appeared quite clear that the class of the differentiable random times should constitute itself an autonomous class possessing the major properties for the purpose of market modeling. This paper is intended to supply a general analysis to this question.

The class is considered of all random times whose conditional distribution functions are differentiable with respect to $\mathbb{F}$ adapted non decreasing processes. Actually this class has been studied in \cite{Libo}. The question raised from \cite{songmodel} compels us to review this study. Accordingly we conclude that the central point for that class of random times is its relation with Cox models, and the results about this class are better presented around this relationship. This idea is just natural (cf. \cite{JS1} for a first example), except Cox model can not exist everywhere (cf. for example \cite{BR, JS2}). It is a situation already encountered in \cite{songthesis} where a systematical use of auxiliary spaces solved the problem. Consequently, the appropriate statement about the relationship between the differentiable random times and the Cox models is that any such random time can be isomorphically implanted into an auxiliary model which is absolutely continuous with respect to a Cox model. 

Here are the mains points of this paper. The method of auxiliary space is essential in this paper. We refer to \cite{songthesis, songlocalsolution} for its general application in the theory of enlargement of filtration. Section \ref{auxi} presents specific properties necessary to make use of this method in this paper. In this same section the notion of i$\!M$ is recalled, which constitutes the very element upon what the whole paper stands. In Section \ref{differential} the differentiability is defined and the first main result is proved, which relates the differentiable random times to the Cox models. A major difference here with respect to the definition in \cite{Libo} is that a certain martingale property was included as part of the definition of the differentiability in \cite{Libo}. For the purpose of applications, we do not assume this martingale property in the definition. Instead, we will establish this martingale property as the consequence of the differentiability. By means of the relationship between the differentiable random times and the Cox models, the easiest way to construct differentiable random times is to make probability changes on Cox models. But such constructions are not always helpful in the practice, because of the lack of means of computations (cf. \cite{w}). In Section \ref{sdetime} we present the results in \cite{songmodel} about the constructions of differentiable random times via stochastic differential equation. We recall that this construction method provides, besides the differentiable models, various models such as Cox models, density hypothesis model (cf. \cite{ejj}), pseudo stopping time model (cf. \cite{NY}), etc. Section \ref{copula} gives a study on the order statistics of differentiable random times via copulas. Starting from Section \ref{cex} the paper lists some main consequences of the differentiability. Three formulas are established: the conditional expectation formula, the optional splitting formula, and the enlargement of filtration formula. It is to note that the use of Cox model in the auxiliary space enable us to have a quick and transparent proofs of these formulas. At last we want to underline that behind the seemingly habitual computations on a market model, there may be hidden measurability problems as indicated in \cite{songsplitting}. This paper tries to define precisely every items appearing in the formulas.

\

\section{Preliminary}

This section gathers some general results useful in the following sections. 

\subsection{Some properties on $\sigma$-algebras}

We recall some facts on the $\sigma$-algebras. Let $E$ a set and $\mathcal{T}$ a $\sigma$-algebra on $E$. For $A\subset E$, we define$$
A\cap\mathcal{T}
=\{B\subset A: \exists C\in\mathcal{T}, B=A\cap C\}.
$$
The family $A\cap\mathcal{T}$ can be used as a family of subsets in $E$, or as a $\sigma$-algebra on $A$. Notice that, if $f$ and $g$ are two maps on $E$ such that $f=g$ on $A$, we have $A\cap\sigma(f)=A\cap\sigma(g)$. We need another fact on the $\sigma$-algebras. Let $\mathcal{T}_n, n\in\mathbb{N}^*$, be a decreasing sequence of $\sigma$-algebras on $E$. Let $F$ be another space and $\eta$ be a map from $F$ into $E$. We have$$
\eta^{-1}(\cap_{n\in\mathbb{N}^*}\mathcal{T}_n)
=\cap_{n\in\mathbb{N}^*}\eta^{-1}(\mathcal{T}_n).
$$ 
Obviously the right hand side term contains the left hand side term. Let $B$ be an element in $\cap_{n\in\mathbb{N}^*}\eta^{-1}(\mathcal{T}_n)$. For any $n\in\mathbb{N}^*$, there exists a $C_n\in \mathcal{T}_n$ such that $B=\eta^{-1}(C_n)$. Let $C=\limsup_{n\rightarrow \infty}C_n\in\cap_{n\in\mathbb{N}^*}\mathcal{T}_n$. We check that $$
\eta^{-1}(C)=\eta^{-1}(\cap_{n\geq 1}\cup_{m\geq n}C_m)
=\cap_{n\geq 1}\cup_{m\geq n}\eta^{-1}(C_m)=B.
$$
This proves the above identity. Consider a subset $D$ of $F$. Applying the previous result with the identity map from $D$ into $F$, we can state the above identity in a general form:

\bl\label{intersection}
We have the identity$$
D\cap\eta^{-1}(\cap_{n\in\mathbb{N}^*}\mathcal{T}_n)
=D\cap(\cap_{n\in\mathbb{N}^*}\eta^{-1}(\mathcal{T}_n))
=\cap_{n\in\mathbb{N}^*}(D\cap\eta^{-1}(\mathcal{T}_n)).
$$ 
\el

\subsection{Negligible sets}

We need the following lemma which describes the completion of a $\sigma$-algebra in term of the $\sigma$-algebra itself. This description will be useful when we compare the completion of a $\sigma$-algebra on the original space with the completion of a $\sigma$-algebra on the auxiliary space.

\bl\label{negtribu}
Let $\mathcal{T}_1, \mathcal{T}_2$ be two $\sigma$-algebras on some common space $\Omega$. Let $\nu$ be a probability measure defined on the two $\sigma$-algebras. Let $\mathcal{N}$ be the family of the $(\nu,\mathcal{T}_2)$ negligible sets. Then,$$
\mathcal{T}_1\vee\sigma(\mathcal{N})
=\{X\subset \Omega: \exists B\in\mathcal{T}_1, A\in\mathcal{T}_2, \nu[A]=1, X\cap A=B\cap A\}.
$$
\el

\textbf{Proof.} Denote the right hand side term of the above formula by $\mathcal{J}$. Then, $\Omega\in\mathcal{J}$. If $X\in\mathcal{J}$, let $B\in\mathcal{T}_1$ and $A\in\mathcal{T}_2$ such that $\nu[A]=1$ and $X\cap A=B\cap A$. Then, $X^c\cap A=B^c\cap A$, which means $X^c\in\mathcal{J}$. If $X_n\in\mathcal{J}$ for $n\in\mathbb{N}^*$, let $B_n\in\mathcal{T}_1$ and $A_n\in\mathcal{T}_2$ such that $\nu[A_n]=1$ and $X\cap A_n=B_n\cap A_n$. Set $A=\cap_{n\in\mathbb{N}^*}A_n$. Then, $\nu[A]=1-\nu[A^c]=1,$ while $$
(\cup_{n\in\mathbb{N}^*}X_n)\cap A
=\cup_{n\in\mathbb{N}^*}(X_n\cap A)
=\cup_{n\in\mathbb{N}^*}(B_n\cap A)
=(\cup_{n\in\mathbb{N}^*}B_n)\cap A),
$$
i.e. $\cup_{n\in\mathbb{N}^*}X_n\in\mathcal{J}$. The family $\mathcal{J}$ is a $\sigma$-algebra.

The $\sigma$-algebra $\mathcal{J}$ contains clearly $\mathcal{T}_1$. It also contains $\mathcal{N}$. Actually, for any $X\in\mathcal{N}$, there exists a $C\in\mathcal{T}_2$ such that $\nu[C]=0$ and $X\subset C$. Let $A=C^c$, we have $\nu[A]=1$ and $X\cap A=\emptyset=\emptyset\cap A$. This means that $X\in\mathcal{J}$. 

On the other hand, for any $X\in\mathcal{J}$, let $B\in\mathcal{T}_1$ and $A\in\mathcal{T}_2$ such that $\nu[A]=1$ and $X\cap A=B\cap A$. Then,$$
X=X\cap A+X\cap A^c
=B\cap A+X\cap A^c\in\mathcal{T}_1\vee\sigma(\mathcal{N}).
$$
This means that $\mathcal{J}\subset \mathcal{T}_1\vee\sigma(\mathcal{N})$. \ok

\subsection{The first zero of a non negative supermartingale}

We need a lemma on the non negative super martingales. 

\bl\label{YR0}
Let $Y$ be a non negative supermartingale defined on some probability space with filtration. Let $X$ be a bounded martingale. Consider the predictable bracket $\cro{X,Y}$. Then, $$
\int_0^\infty \ind_{\{Y_{s-}=0\}}d\cro{X,Y}_s=0.
$$
Let $V$ be the drift part of $Y$ in its canonical decomposition ($V$ being decreasing). Let $R=\inf\{t\in\mathbb{R}_+:\ Y_t=0\}$ and $R'=R\ind_{\{Y_{R-}=0\}}+\infty\ind_{\{Y_{R-}>0\}}$. Then, $\Delta_{R'}V=0$ on $\{R'<\infty\}$.  
\el

\textbf{Proof.} Let $\mathtt{F}=\inf\{0<R<\infty, Y_{R-}=0\}$. By \cite[Theorem 2.62, Corollary 12.5]{Yan} (also cf.\cite[Theorem 3.35]{Yan}), $$
\{Y_-=0\}=\mathtt{F}\cap[R,\infty)+\mathtt{F}^c\cap(R,\infty).
$$
Hence, $$
[R_\mathtt{F}]=[0,R]\cap \{Y_-=0\}
$$
is a predictable set, i.e. $R_\mathtt{F}$ is a predictable stopping time. We compute now
$$
\dcb
\int_0^\infty \ind_{\{Y_{s-}=0\}}d\cro{X,Y}_s
&=\int_0^\infty \ind_{\{0<s\leq R\}}\ind_{\{Y_{s-}=0\}}d\cro{X,Y}_s	\\
&=\Delta_{R_{\mathtt{F}}}\cro{X,Y}\
=\mathbb{E}[\Delta_{R_{\mathtt{F}}}X \ \Delta_{R_{\mathtt{F}}}Y|\mathcal{F}_{R_{\mathtt{F}}}]\
=0.
\dce
$$
This proves the first part of the lemma. 

As for the second part of the lemma, we note that $Y_-={^{p}(Y)}-\Delta V$ (where the superscript $^p$ denote the optional projection). So, $Y_{s-}>0$ for a $0\leq s<\infty$, whenever $\Delta_sV<0$. However, $Y_{R'-}=0$ on $\{R'<\infty\}$ so that $\Delta_{R'}V=0$ on $\{R'<\infty\}$. \ok

\

\subsection{Change of variable}\label{changevariable}

Let $a$ be a real non negative right continuous non deceasing function on $\mathbb{R}_+$. Following \cite{Yan} we introduce the right-inverse of $a$:$$
c(s)=\inf\{u\in\mathbb{R}_+: a(u)>s\},\ s\in\mathbb{R}_+.
$$
Then, $c$ is a non negative right continuous non decreasing function. We have$$
\dcb
\{a(u)\geq s\}=\{c(s-)\leq u\},\\
\{a(u-)\leq s\}=\{c(s)\geq u\},\\
a(c(s)-)\leq s\leq a(c(s-)),\\
a(u)=\inf\{s\in\mathbb{R}_+: c(s)>u\}.
\dce
$$
For any non negative Borel function $f$, we have the identity:$$
\int_{[0,\infty)}f(u)da(u)
=\int_{[0,\infty)}\ind_{\{c(s-)<\infty\}}f(c(s-))ds.
$$ 
In particular, for any $t\in\mathbb{R}_+$,$$
\dcb
\int_{[0,t]}e^{-a(u)}da(u)
=\int_{[0,\infty)}\ind_{\{c(s-)\leq t\}}e^{-a(c(s-))}ds
\leq\int_{[0,\infty)}\ind_{\{s\leq a(t)\}}e^{-s}ds=1-e^{-a(t)}< 1.
\dce
$$

\

\section{Random time and the progressive enlargement of filtration}\label{auxi}

The basic setting of this work is a stochastic structure $(\Omega,\mathcal{A},\mathbb{F},\mathbb{P})$ equipped with a $\tau$, where $(\Omega,\mathcal{A})$ is a measurable space, and $\mathbb{P}$ is a probability measure on $\mathcal{A}$, and $\mathbb{F}$ is a right continuous filtration $\mathbb{F}=(\mathcal{F}_t)_{t\in\mathbb{R}_+}$ (with $\mathcal{F}_\infty=\vee_{t\in\mathbb{R}_+}\mathcal{F}_t\subset \mathcal{A}$) which contain the $(\mathbb{Q},\mathcal{F}_\infty)$ null sets, and finally $\tau$ is a random variable taking values in $[0,\infty]$ (that will be called random time below). Let us denote by $\nmid$ the map defined on $[0,\infty]^2$ into $[0,\infty]$ such that, for $a,b\in[0,\infty]$, $a\nmid b=a$ if $a\leq b$ and $a\nmid b=\infty$ if $a> b$. 
We introduce then $\mathbb{G}^0=(\mathcal{G}^0_t)_{t\in\mathbb{R}_+}$ the filtration defined by $\mathcal{G}^0_t=
\cap_{s>t}(\mathcal{F}_s\vee\sigma(\tau\nmid s))
$,
and its completion $\mathbb{G}=(\mathcal{G}_t)_{t\in\mathbb{R}_+}$ where the $\sigma$-algebra $\mathcal{G}_t$ is the $\sigma$-algebra $\mathcal{G}^0_t$ completed with the $(\mathbb{Q},\mathcal{F}_\infty\vee\sigma(\tau))$ null sets. The filtration  $\mathbb{G}$ is the progressive enlargement of the filtration $\mathbb{F}$ with the random time $\tau$.

\

\subsection{Product measurable space}

The fundamental idea in \cite{songthesis} (see also \cite{songlocalsolution} for a recent presentation) to deal with the enlargement of filtrations is the following two steps scheme, called algorithm in that thesis : firstly, to introduce another probability measure $\mathbb{Q}'$ on $\mathcal{G}_\infty$ under which the random time is sufficiently independent of the filtration $\mathbb{F}$ so that the problems in $\mathbb{G}$ has a solution; and then, to link the $(\mathbb{Q},\mathbb{G})$ semimartingales with the $(\mathbb{Q}',\mathbb{G})$ semimartingales through Girsanov's theorem. But in fact, generally this two steps scheme will not be applied directly on the original probability space, because, for a general probability space, it may impossible to define a new probability measure $\mathbb{Q}'$ on $\mathcal{G}_\infty$ which makes $\tau$ independent of $\mathbb{F}$. Instead there is a third step in this methodology, i.e. to introduce an auxiliary space on which we apply the two steps scheme, and to apply the "invariance principle" to go back to the original space. The "invariance principle" means that the enlargement of filtration problem is a problem in "law". The problem has a solution whenever it has a solution on an isomorphic (in large sense) auxiliary space.

Concretely for the progressive enlargement of filtration, we consider the map $\phi(\omega)= (\omega,\tau(\omega))$ from $\Omega$ into the product space $\Omega\times[0,\infty]$ and we define the auxiliary space as the product space $\Omega\times[0,\infty]$ equipped with the product $\sigma$-algebra
$\mathcal{A}\otimes\mathcal{B}[0,\infty]$ and with the image probability measure : $\check{\mathbb{Q}}[B]=\mathbb{Q}[\phi^{-1}(B)]$, $B\in\mathcal{A}\otimes\mathcal{B}[0,\infty]$. To have a representation of the probability structure $(\mathbb{Q},\mathbb{F})$ on the product space, we consider the maps $\pi(\omega,u)=\omega$, $(\omega,u)\in\Omega\times[0,\infty]$. With the map $\pi$ we draw the filtration $\mathbb{F}$ onto the product space $\Omega\times[0,\infty]$ as follows:$$
\check{\mathbb{F}}=(\check{\mathcal{F}}_t:t\geq 0)=\pi^{-1}(\mathbb{F})=(\pi^{-1}(\mathcal{F}_t):t\geq 0).
$$
To have a representation of the random time $\tau$, we introduce the map $\hat{\tau}(\omega,u)=u$. We check that, for any $B\in\mathcal{F}_\infty$, for $a\in[0,\infty]$, we have $\pi^{-1}(B)\in\check{\mathcal{F}}_\infty$ and $$
\mathbb{Q}[B\cap\{\tau\leq a\}]
=\check{\mathbb{Q}}[\pi^{-1}(B)\cap\{\hat{\tau}\leq a\}].
$$
With the filtration $\check{\mathbb{F}}$ and the random time $\hat{\tau}$ we define on the product space the primal enlargement and the progressive enlargement, i.e. $\check{\mathbb{G}}^0=(\check{\mathcal{G}}^0_t)_{t\geq 0}$ defined by $\check{\mathcal{G}}^0_t
=\cap_{s>t}(\check{\mathcal{F}}_s\vee\sigma(\hat{\tau}\nmid s))$, and $\check{\mathbb{G}}=(\check{\mathcal{G}}_t)_{t\geq 0}$ where $\check{\mathcal{G}}_t$ is the completion of $\check{\mathcal{G}}^0_t$ by the $\check{\mathbb{Q}}$ null sets with respect to $\check{\mathcal{F}}_\infty\vee\sigma(\hat{\tau})$. Now we have a representative $(\check{\mathbb{Q}}, \check{\mathbb{F}},\check{\mathbb{G}},\check{\tau})$ on the product space which duplicates the original quadruplet $(\mathbb{Q},\mathbb{F},\mathbb{G},\tau)$. For any map $X$ defined on the space $\Omega$, we define the map $\check{X}=X(\pi)$ on the product space.

\brem
We note that there exist two copies on the product space which duplicate the random time $\tau$, i.e. the maps $\check{\tau}$ and $\hat{\tau}$. These two representatives are indistinguishable under the probability $\check{\mathbb{Q}}$. However, when one wants to construct a new probability in introducing some independence into the progressive enlargement of filtration, we deal only with the representative $\hat{\tau}$.  
\erem

\brem
Notice that in stochastic calculus, one usually omits to write the dependence on $\omega$ of random variables. We adopt this convention. Therefore in this paper, a stochastic process $X(t,\omega)$ is written as $X_t$, and a function $f(\omega,u)$ on the product space is written as $f(u)$, so on.
\erem

\

\subsection{Measurability relationship between the original and the auxiliary spaces}

Now we compare the different optional $\sigma$-algebra on the original and on the auxiliary spaces.

\bl\label{optionalinclusion}
Let $\Phi$ to be the map from $\mathbb{R}_+\times\Omega$ into $\mathbb{R}_+\times(\Omega\times[0,\infty])$ with $\Phi(t,\omega)=(t,\phi(\omega))$. We have the following inclusion relationship:$$
\Phi^{-1}(\mathcal{O}(\check{\mathbb{G}}^0))
\subset \mathcal{O}(\mathbb{G}^0)\subset
\Phi^{-1}(\mathcal{O}(\check{\mathbb{G}}))
\subset \mathcal{O}(\mathbb{G}).
$$
\el

\textbf{Proof.}
The following relationships hold:$$
\phi^{-1}\left(\check{\mathcal{F}}_t\vee\sigma(\hat{\tau}\nmid t)\right)
=
\mathcal{F}_t\vee\sigma(\tau\nmid t),\ \
\phi^{-1}(\check{\mathcal{G}}^0_t)
=
\mathcal{G}^0_t,
$$
for every $t\in\mathbb{R}_+$, and hence by Lemma \ref{negtribu},
$
\phi^{-1}(\check{\mathcal{G}}_t)
\subset
\mathcal{G}_t.
$
Then, for any (everywhere) càdlàg $\check{\mathbb{G}}^0$ (respectively $\check{\mathbb{G}}$) adapted process $X$, $X\circ\Phi$ is a càdlàg $\mathbb{G}^0$ (respectively $\mathbb{G}$) adapted process. We have therefore$$
\Phi^{-1}(\mathcal{O}(\check{\mathbb{G}}^0))
\subset \mathcal{O}(\mathbb{G}^0),\
\Phi^{-1}(\mathcal{O}(\check{\mathbb{G}}))
\subset \mathcal{O}(\mathbb{G}).
$$
For a (everywhere) càdlàg $\mathbb{G}^0$ adapted process $X$ bounded by 0 and 1, for any $t\in\mathbb{Q}_+$, let $f_t$ be a $\check{\mathcal{G}}^0_t$ measurable function bounded by 0 and 1, such that $X_t=f_t(\phi)$ (cf. \cite[Theorem 1.5]{Yan}). 
For every $a\in\mathbb{R}_+$, let $\mathtt{W}'_a$ to be the set of $(\omega,u)$ in the product space such that, there exists a $\epsilon=\epsilon(\omega,u)>0$, and $(f_t(\omega,u): t\in\mathbb{Q}_+\cap[0,a+\epsilon))$ is the restriction on $\mathbb{Q}_+\cap[0,a+\epsilon)$ of a càdlàg function. Then, according to \cite[Chapitre IV $n^\circ$18]{DM}, $\mathtt{W}'_a$ is in $\cap_{\epsilon>0}\check{\mathcal{G}}^0_{a+\epsilon}=\check{\mathcal{G}}^0_a$. Since $X$ is càdlàg, $\check{\mathbb{Q}}[\mathtt{W}'_a]=1$. Notice that $\mathtt{W}'_a\supset \mathtt{W}'_b$ for any $a\leq b$ and $\cup_{a<b<T}\mathtt{W}'_b=\mathtt{W}'_a$. We define$$
\hat{X}_t=\liminf_{s\rightarrow t, s>t, s\in\mathbb{Q}_+}f_s\ind_{\mathtt{W}'_s},\ \ t\in\mathbb{R}_+.
$$
There are two cases. Firstly, if $(\omega,u)$ belongs to no of the $\mathtt{W}'_s$ for $s>t$, we have $\hat{X}_v(\omega,u)=0$ for all $v\geq t$. Hence, $\hat{X}(\omega,u)$ is right continuous at $t$. Secondly, if $(\omega,u)$ belongs to some $\mathtt{W}'_{s_0}$ for a $s_0>t$, then, $(\omega,u)$ belongs to all $\mathtt{W}'_s$ for a $t<s\leq s_0$, and $f(\omega,u)$ on $\mathbb{Q}_+\cap[0,s_0)$ is the restriction of a càdlàg function $g$ on $\mathbb{Q}_+\cap[0,s_0)$. This implies that $\hat{X}_v(\omega,u)=g(v)$ for $0\leq v<s_0$. In particular, $\hat{X}(\omega,u)$ is right continuous at $t$. We just prove that $\hat{X}$ is a right continuous process. Moreover, we check directly that $\phi(\omega)\in\mathtt{W}'_s$ for all $s\geq 0,\omega\in\Omega$, and therefore $X=\hat{X}(\phi)$. Notice that $\hat{X}\in\mathcal{O}(\check{\mathbb{G}})$ (cf. \cite[Theorem 4.32]{Yan}). This being true for an (everywhere) càdlàg $\mathbb{G}^0$ adapted process $X$ bounded by 0 and 1, we prove the last inclusion relation which was missing in the above inclusion sequence. \ok

\bl\label{oopp}
Let $a\in\mathbb{R}_+$. Let $H$ be a non negative $\check{\mathbb{G}}$ optional process. Let $A$ be a non decreasing $\mathbb{F}$ optional process. Consider $H$ as a function of three variables $H_t(\omega,u)$, $t\in\mathbb{R}_+$, $(\omega,u)\in\Omega\times[0,\infty]$, and denote by $H_t(u)$ the map $\omega\rightarrow H_t(\omega,u)$. Then, the process $\int_{[0,a]}H_t(v)dA_v$, $t\in[a,\infty)$, is $\mathbb{F}$ optional on the interval $[a,\infty)$. 
\el

\textbf{Proof.} Consider $\mathbf{H}$ the set of bounded $\check{\mathbb{G}}$ optional processes $H$ such that the statement of the lemma is valid on $H$. It is clear that $\mathbf{H}$ is a functional monotone class in the sense of \cite[Theorem (3.2)]{RW}. $\mathbf{H}$ is therefore a vector space containing the constant functions and closed under uniform convergence. Let $\mathbf{H}^0$ be the set of all bounded $\check{\mathbb{G}}$ adapted càdlàg processes. $\mathbf{H}^0$ is closed under multiplication. Let us show that any element $H$ in $\mathbf{H}^0$ is an element in $\mathbf{H}$. Actually, for fixed $t\in\mathbb{R}_+$, by an argument by monotone class theorem, the random variable $\int_{[0,a]}H_t(v)dA_v$ is $\mathcal{F}_{t+\epsilon}$ measurable for any $\epsilon>0$, hence it is $\mathcal{F}_{t}$ measurable. In addition, the process $\int_{[0,a]}H_t(v)dA_v, t\in\mathbb{R}_+$, is càdlàg by dominated convergence theorem.

By $\cite{RW}$ $\mathbf{H}$ contains all bounded $\check{\mathbb{G}}$ optional process. For a general non negative $\check{\mathbb{G}}$ optional process, it is the increasing limit of a sequence of bounded $\check{\mathbb{G}}$ optional process. This proves the lemma. \ok

\

\subsection{i$\! M$ increasing family of martingales}

In mathematical modeling of financial market through the progressive enlargement of filtration, the most important characteristic of the random time $\tau$ should be its conditional distribution function. This notion is formalized in \cite{JS2} with the following definition.

\bd
An i$\!M(\mathbb{Q},\mathbb{F})$ family (or simply i$\!M$ family) is a family of processes $(M^u: u\in[0,\infty])$ satisfying the following conditions: 

\ebe 
\item[$1\cdot$] 
For every $u\in[0,\infty]$, $M^u$ is a $(\mathbb{Q},\mathbb{F})$ martingale on
$[u,\infty]$ taking values in $[0,1]$. 
\item[$2\cdot$] 
For every $t\in[0,\infty]$, the random map $u\in[0,t]\rightarrow M^u_t$ is a right continuous non-decreasing function. 
\item[$3\cdot$]
$M^\infty_\infty=1$.
\dbe
If the condition $1\cdot$ and $2\cdot$ are replaced by
\ebe
\item[$1'\cdot$]
For every $u\in[0,\infty]$, $M^u$ is a $(\mathbb{Q},\mathbb{F})$ martingale on
$[0,\infty]$ taking values in $[0,1]$.
\item[$2'\cdot$] 
For every $t\in[0,\infty]$, the random map $u\in[0,\infty]\rightarrow M^u_t$ is a right continuous non-decreasing function.
\dbe
we say that the i$\!M$ is complete.
\ed

The following theorem is borrowed from \cite{JS2}.

\bethe \label{theoremIM}  
\ebe 
\item[$i\cdot$]
For any random time $\tau$ on the filtrationed probability space $(\Omega,\mathbb{F},\mathbb{Q})$, there exists a complete i$\! M$ family, denoted by $(M^u:u\in[0,\infty])$, such that, for $u,t\in[0,\infty]$, $$
M^u_t=\mathbb{Q}[\tau\leq u|\mathcal{F}_t].
$$
This family is unique in the sense that, if $(\widetilde{M}^u:u\in[0,\infty])$ is another i$\! M$ family satisfying the above condition, there exists a $\mathbb{Q}$ null set $A$ such that, for $\omega\notin A$, $M^u_t(\omega)=\widetilde{M}^u_t(\omega)$ for all $0\leq u, t\leq \infty$. We say that $(M^u:u\in[0,\infty])$ is associated with triplet $(\mathbb{Q},\mathbb{F},\tau)$.

\item$ii\cdot$
Let $(M^u: u\in[0,\infty])$ be an i$M(\mathbb{Q},\mathbb{F})$. There is a
unique probability measure $\mathbb{Q}'$ on
the product measurable space, which coincides with $\check{\mathbb{Q}}$ on $\check{\mathcal{A}}=\pi^{-1}(\mathcal{A})$ and satisfies $\mathbb{Q}'[\check{\tau}\leq u|\check{\mathcal{F}}_t]=M^u_t(\pi)$ for $0\leq
u\leq t\leq  \infty$, and $\mathbb{Q}'[\check{\tau}\leq u|\check{\mathcal{A}}]=M^u_\infty(\pi)$. We call $\mathbb{Q}'$ the measure on the product space
associated with the triplet $(\mathbb{Q},\mathbb{F},i\!M)$.

\item$iii\cdot$
As a consequence of $ii\cdot$, for any i$\!M$ family, it has a compete extension. 
\dbe
\ethe

We also need the following technical lemmas to deal with the i$\! M$ family. The first lemma is a direct consequence of the definition and of the monotone class theorem.

\bl\label{Mint}
Let $(M^u:u\in[0,\infty])$ be the i$M$ family associated with the triplet $(\mathbb{Q},\mathbb{F},\tau)$. For any $t\in[0,\infty]$, for any non negative $\mathcal{F}_t\otimes \mathcal{B}[0,\infty]$ function $f$, we have$$
\mathbb{E}[f(\tau)]=\mathbb{E}[\int_{[0,\infty]}f(u)d_uM^u_t].
$$
\el

\bl
For any i$\! M(\mathbb{Q},\mathbb{F})$ family $(M^u:u\in[0,\infty])$, the maps $(M^t_t)_{t\in\mathbb{R}_+}$ defines a $\mathbb{F}$ optional process.
\el

\textbf{Proof.} We take the complete extension of the i$\! M$ family $(M^u:u\in[0,\infty])$. Then the map $$
((t,\omega),u)\in(\mathbb{R}_+\times\Omega)\times[0,\infty]\rightarrow M^u_t(\omega)
$$ 
is $\mathcal{O}(\mathbb{F})\times\mathcal{B}[0,\infty]$ measurable, because of the right continuity in $u$. Define a map $\psi$ from $\mathbb{R}_+\times\Omega$ into $(\mathbb{R}_+\times\Omega)\times[0,\infty]$ by $\psi(t,\omega)=((t,\omega), t)$. Then, $$
\psi^{-1}(\mathcal{O}(\mathbb{F})\times\mathcal{B}[0,\infty])
=\mathcal{O}(\mathbb{F}).
$$
The lemma comes as a consequence. \ok

\bl\label{vaetvien}
Let $\mathbf{M}=(M^u:u\in[0,\infty])$ be the i$\! M$ family associated with $(\mathbb{Q},\mathbb{F},\tau)$. Let $\widetilde{\mathbf{M}}=(\hat{M}^u:u\in[0,\infty])$ be the i$\! M$ family on the product space associated with $(\check{\mathbb{Q}},\check{\mathbb{F}},\hat{\tau})$. Then, $\mathbf{M}(\pi)$ is a version of $\widetilde{\mathbf{M}}$, and $\widetilde{\mathbf{M}}(\phi)$ is a version of $\mathbf{M}$.
\el

\textbf{Proof.} The lemma is the consequence of the right continuity of $\widetilde{\mathbf{M}}$ and $\mathbf{M}$, and of the following identity. For $u,t\in[0,\infty]$ with $u\leq t$, for any $B\in\mathcal{F}_t$, we have$$
\check{\mathbb{E}}[\ind_B(\pi)M^u_t(\pi)]
=\mathbb{E}[\ind_B M^u_t]
=\mathbb{E}[\ind_B\ind_{\{\tau\leq u\}}]
=\check{\mathbb{E}}[\ind_B(\pi)\ind_{\{\hat{\tau}\leq u\}}]
=\check{\mathbb{E}}[\ind_B(\pi)\widetilde{M}^u_t]
=\mathbb{E}[\ind_B\widetilde{M}^u_t(\phi)].\ \ok
$$

\

\section{Differentiable i$\! M$ and Cox measure}\label{differential}

From now on, we fix a constant $T\in[0,\infty]$. Here is the definition of the differentiability of an i$\! M$.

\bd\label{ac}
Let $(M^u:u\in[0,\infty])$ be an i$\! M(\mathbb{Q},\mathbb{F})$ family.	Let $A$ be a non negative $\mathbb{F}$ adapted increasing càdlàg process. For $t\in[0,\infty]$, $(M^u:u\in[0,\infty])$ is said to be differentiable at $t$ with respect to $A$, if there exists a non negative $\mathcal{F}_t\otimes\mathcal{B}[0,\infty]$ measurable function $\mathsf{p}_t(\omega,v)$ such that, for almost all $\omega$, $$
M^u_t=\int_{[0,u]}\mathsf{p}_t(v)dA_v,\ \forall u\in[0,t].
$$
(As usual we omit $\omega$.) We call $\mathsf{p}_t$ a density function at $t$. If $(M^u:u\in[0,\infty])$ is differentiable at every $t\in[0,T)$ with respect to the same increasing process $A$, we say that $(M^u:u\in[0,\infty])$ is differentiable on $[0,T)$. 
\ed

\brem\label{aaa}
If we replace $dA_v$ by $e^{-A_v}dA_v$ and $\mathsf{p}_t(v)$ by $\mathsf{p}_t(v)e^{A_v}$, we can assume in the above definition that $A_{t}< 1$ for $t\in\mathbb{R}_+$ (cf. subsection \ref{changevariable}). 
\erem

The following lemma is the consequence of Lemma \ref{vaetvien}. 

\bl\label{passpass}
Let $\mathbf{M}=(M^u:u\in[0,\infty])$ be the i$\! M$ family associated with $(\mathbb{Q},\mathbb{F},\tau)$. Let $\widetilde{\mathbf{M}}=(\hat{M}^u:u\in[0,\infty])$ be the i$\! M$ family on the product space associated with $(\check{\mathbb{Q}},\check{\mathbb{F}},\hat{\tau})$. Let $A$ be a non negative $\mathbb{F}$ adapted increasing process. Define $\check{A}=A(\pi)$. 

If $\mathbf{M}$ is differentiable with respect to $A$ on $[0,T)$ with a density function $\mathsf{p}$, then $\widetilde{\mathbf{M}}$ is differentiable with respect to $\check{A}$ on $[0,T)$ with the density function $\check{\mathsf{p}}$ defined by
$$
\check{\mathsf{p}}_t((\omega,u),v)=\mathsf{p}_t(\pi(\omega,u),v),\ (\omega,u)\in\Omega\times[0,\infty], v\in[0,\infty].
$$
Conversely, if $\widetilde{\mathbf{M}}$ is differentiable with respect to $\check{A}$ on $[0,T)$ with a density function $\check{\mathsf{p}}$, then $\mathbf{M}$ is differentiable with respect to $A$ on $[0,T)$ with the density function $\mathsf{p}$ defined by
$$
\mathsf{p}_t(\omega,v)=\check{\mathsf{p}}_t(\phi(\omega),v),\ \omega\in\Omega, v\in[0,\infty].
$$
\el

We are going to display the various formulas in term of the density function $\mathsf{p}$. However, to really be able to do so, we need first of all a modified version of $\mathsf{p}$. We introduce an additional notion (cf. \cite[Chapter 3]{lisbonn} for the notion of Cox process). 

\bd
We call a probability measure $\mathbb{Q}^0$ on $\mathcal{F}_\infty\vee\sigma(\tau)$ a Cox measure with respect to $\mathbb{Q}|_{\mathcal{F}_\infty}$ (the restriction of $\mathbb{Q}$ on $\mathcal{F}_\infty$), if there exists a non negative $\mathbb{F}$ adapted increasing càdlàg process $A$ such that 
\ebe
\item[i.]
the two probability measures $\mathbb{Q}^0=\mathbb{Q}$ on $\mathcal{F}_\infty$;
\item[ii.]
for all $u,t\in[0,\infty]$ with $u\leq t$, we have $\mathbb{Q}^0[\tau\leq u|\mathcal{F}_t]=A_u$. 
\dbe
We also say that $\mathbb{Q}^0$ is the Cox measure associated with $(\mathbb{Q}|_{\mathcal{F}_\infty},\mathbb{F},A,\tau)$.
\ed

\brem\label{coxim}
Note that the i$\! M$ family associated with $(\mathbb{Q}^0,\mathbb{F},\tau)$ is given by $M^u_t=A_u$ for $0\leq u\leq t\leq \infty$.
\erem

\bl
For any non negative $\mathbb{F}$ adapted non decreasing càdlàg process $A=(A_t,t\in\mathbb{R}_+)$ such that $A_{\infty-}\leq 1$, we extend $A$ to the domain $[0,\infty]$ by defining $A_\infty=1$. Let $\check{A}=A(\pi)$. Then, there exists a Cox measure $\check{\mathbb{Q}}^0$ associated with $(\check{\mathbb{Q}}|_{\check{\mathcal{F}}_\infty},\check{\mathbb{F}},\check{A},\hat{\tau})$ on the product space.
\el

\textbf{Proof.} It is necessary and sufficient to define the Cox measure as follows:
$$
\check{\mathbb{Q}}^0[h]
=\int_{w\in\Omega\times[0,\infty]} d\mathbb{Q}(\omega)\int_{[0,\infty]}h(\omega,v)d\check{A}_v(\pi(w)),
$$
for non negative $\mathcal{F}_\infty\otimes\mathcal{B}[0,\infty]$ measurable function $h$. \ok

The following theorem, in which a very precise version of the density function $\mathsf{p}$ is studied, is essential for this paper. We consider this theorem as a specific version of the Follmer's lemma (cf. \cite[Theorem 2.44]{Yan}) adapted to the case of a progressively enlarged filtration.

\bethe\label{acp}
Let $Z$ be the $(\mathbb{Q},\mathbb{F})$ supermartingale $\mathbb{Q}[t<\tau|\mathcal{F}_t], t\in\mathbb{R}_+$ (called the Azéma supermartingale of $\tau$). Let $A$ be a non negative $\mathbb{F}$ adapted increasing càdlàg process such that $A_\infty=1$ and $A_{t}< 1$ for $t\in\mathbb{R}_+$. Let $(M^u: u\in[0,\infty])$ be an i$\! M(\mathbb{Q},\mathbb{F})$ family which is differentiable on $[0,T)$ with respect to $A$ with a density function $\mathsf{p}$. Let $\check{\mathbb{Q}}^0$ be the Cox measure on the product space associated with $(\check{\mathbb{Q}}|_{\check{\mathcal{F}}_\infty},\check{\mathbb{F}},\check{A},\hat{\tau})$.

Then, there exists a three variable function $\mathsf{p}_{t+}(\omega,u), (\omega,u)\in\Omega\times[0,\infty], t\in[0,T)$ such that 
\ebe
\item
the process $\mathsf{p}_{t+}$, $t\in[0,T)$, is $\check{\mathbb{G}}^0$ adapted on the interval $[0,T)$, and, for all $(\omega,u)\in\Omega\times[0,\infty]$, the map $t\rightarrow \mathsf{p}_{t+}(\omega,u)$ is everywhere càdlàg on $[0,T)$; 

\item
$\check{\mathbb{Q}}$ is absolutely continuous with respect to $\check{\mathbb{Q}}^0$ on $\check{\mathcal{G}}_t$ for every $t\in[0,T)$ and the process$$
\mathsf{P}_t=\ind_{\{t<u\}}\frac{Z_t(\omega)}{1-A_t(\omega)}+
\ind_{\{u\leq t\}}\mathsf{p}_{t+}(\omega,u),\ t\in[0,T),
$$
is the corresponding density process;

\item
$M^a_t$ and $\int_{[0,a]}\mathsf{p}_{t+}(v)dA_v$ for $0\leq a\leq t<T$ are $\mathbb{Q}$ distinguishable as two two-parameter processes, where $\mathsf{p}_{t+}(v)$ denotes the map $\omega\rightarrow \mathsf{p}_{t+}(\omega,v)$. In particular $t\rightarrow \int_{[0,a]}\mathsf{p}_{t+}(v)dA_v$ is càdlàg on $[a,T)$. 
\dbe
\ethe

\brem
We will keep the notation $\mathsf{p}$ to denote the original version of the density function. We will use $\mathsf{p}_+$ to denote the version of the density function established in this theorem. Later we will have a third version of the density function.
\erem

\textbf{Proof.} We use the notations in Lemma \ref{passpass}. For every $t\in[0,T)$, for any non negative bounded Borel functions $h$ on $[0,\infty]$, for any $B\in\mathcal{F}_t$, applying Lemma \ref{Mint}, we have$$
\dcb
\check{\mathbb{E}}^0[\ind_B(\pi)h(\hat{\tau}\nmid t)\ind_{\{\hat{\tau}\leq t\}}]
&=\mathbb{E}[\ind_Bh(\tau\nmid t)\ind_{\{\tau\leq t\}}]
=\mathbb{E}[\ind_B\int_{[0,t]}h(u)d_uM^u_t]\\

&=\mathbb{E}[\ind_B\int_{[0,t]}h(u)\mathsf{p}_t(u)dA_u]

=\check{\mathbb{E}}^0[\ind_B(\pi)h(\hat{\tau}\nmid t)\ind_{\{\hat{\tau}\leq t\}}\mathsf{p}_t].
\dce
$$
and
$$
\dcb
\check{\mathbb{E}}^0[\ind_B(\pi)h(\hat{\tau}\nmid t)\ind_{\{t<\hat{\tau}\}}]
&=\mathbb{E}[\ind_Bh(\tau\nmid t)\ind_{\{t<\tau\}}]
=\mathbb{E}[\ind_Bh(\infty)\ind_{\{t<\tau\}}]
=\mathbb{E}[\ind_Bh(\infty)Z_t]\\

&=\mathbb{E}[\ind_Bh(\infty)\frac{Z_t}{1-A_t}(1-A_t)]\\

&=\check{\mathbb{E}}^0[\ind_B(\pi)h(\infty)\frac{\check{Z}_t}{1-\check{A}_t}\ind_{\{t<\hat{\tau}\}}]

=\check{\mathbb{E}}^0[\ind_B(\pi)h(\hat{\tau}\nmid t)\frac{\check{Z}_t}{1-\check{A}_t}\ind_{\{t<\hat{\tau}\}}].
\dce
$$
Combining these two identities, we obtain$$
\check{\mathbb{E}}[\ind_B(\pi)h(\hat{\tau}\nmid t)]
=\check{\mathbb{E}}^0[\ind_B(\pi)h(\hat{\tau}\nmid t)\left(
\ind_{\{t<\hat{\tau}\}}\frac{\check{Z}_t}{1-\check{A}_t}+
\ind_{\{\hat{\tau}\leq t\}}\mathsf{p}_t
\right)].
$$
This means that $\check{\mathbb{Q}}$ is absolutely continuous with respect to $\check{\mathbb{Q}}^0$ on $\check{\mathcal{F}}_t\vee\sigma(\hat{\tau}\nmid t)$ with density$$
\ind_{\{t<u\}}\frac{Z_t(\omega)}{1-A_t(\omega)}+
\ind_{\{u\leq t\}}\mathsf{p}_t(\omega,u).
$$
The above process is a $\check{\mathbb{Q}}^0$ martingale in the filtration $\check{\mathcal{F}}_t\vee\sigma(\hat{\tau}\nmid t), t\in[0,T)$ (in the primitive sense of martingale). By martingale convergence theorem (cf. \cite[Theorem 2.44 \textit{Follmer's lemma}]{Yan} or \cite[Régularité des trajectoires p.142]{DM}), there exists an everywhere right continuous $(\check{\mathbb{Q}}^0,\check{\mathbb{G}}^0)$ martingale $\mathsf{P}$ on the time interval $[0,T)$ such that, for $\check{\mathbb{Q}}^0$ almost all $(\omega,u)$, 
$$
\mathsf{P}_t(\omega,u)=\lim_{s\in\mathbb{Q}_+:s\downarrow t}
\left(
\ind_{\{s<u\}}\frac{Z_s(\omega)}{1-A_s(\omega)}+
\ind_{\{u\leq s\}}\mathsf{p}_s(\omega,u)
\right),\ \forall t\in[0,T).
$$
Let $\mathsf{p}_{t+}(\omega,u)=\ind_{\{u\leq t\}}\mathsf{P}_t(\omega,u)$. Then, $\mathsf{p}_{t+}$ is everywhere right continuous in $t\in[0,T)$, and
$$
\mathsf{P}_t(\omega,u)=\ind_{\{t<u\}}\frac{Z_t(\omega)}{1-A_t(\omega)}+
\ind_{\{u\leq t\}}\mathsf{p}_{t+}(\omega,u)
$$
is the density function of $\check{\mathbb{Q}}$ with respect to $\check{\mathbb{Q}}^0$ on $\check{\mathcal{G}}^0_t$.

For $a\in[0,T)$, let $S$ be a $\mathbb{F}$ stopping time $S\in[a,T-\epsilon]$ for some $\epsilon>0$. For $B\in\mathcal{F}_S$,
$$
\dcb
\mathbb{E}[\ind_BM^a_S]
&=\mathbb{E}[\ind_B\ind_{\{\tau\leq a\}}]
=\check{\mathbb{E}}^0[\ind_B(\pi)\ind_{\{\hat{\tau}\leq a\}}\mathsf{p}_{\check{S}+}]
=\mathbb{E}[\ind_B\int_{[0,a]}\mathsf{p}_{S+}(v)dA_v],
\dce
$$
where the last equality comes from the definition of Cox measure $\check{\mathbb{Q}}^0$. The process $(\int_{[0,a]}\mathsf{p}_{t+}(v)dA_v: t\in[a,\infty])$ being $\mathbb{F}$ optional according to Lemma \ref{oopp}, by the section theorem (cf. \cite{Yan}), we conclude that the two processes $M^a$ and $(\int_{[0,a]}\mathsf{p}_{t+}(v)dA_v: t\in[a,\infty])$ are $\mathbb{Q}$ indistinguishable on $[a,T)$. In particular $t\rightarrow \int_{[0,a]}\mathsf{p}_{t+}(v)dA_v$ is càdlàg on $[a,T)$. By the right continuity in the variable $a$, $M^a_t$ and $\int_{[0,a]}\mathsf{p}_{t+}(v)dA_v$ for $0\leq a\leq t<T$ also are distinguishable as two two-parameter processes.

Let us now find a version of $\mathsf{p}_+$ which has everywhere left limit on $[0,T)$. For $b\in[0,T)$ let $\mathtt{W}'_b$ (cf. the proof of Lemma \ref{optionalinclusion}) be the set of $(\omega,u)$ such that, there exists a $0<\epsilon=\epsilon(\omega,u)<T-b$, and the restriction of the process $\mathsf{p}_{t+}(\omega,u)$ on $t\in\mathbb{Q}_+\cap[0,b+\epsilon]$ is equal to the restriction on $\mathbb{Q}_+\cap[0,b+\epsilon]$ of a càdlàg function on $[0,T)$. According to \cite[Chapitre IV, n$^\circ$18]{DM} the set $\mathtt{W}'_b$ is in $\cap_{\epsilon>0}\check{\mathcal{G}}^0_{b+\epsilon}=\check{\mathcal{G}}^0_b$ and $\check{\mathbb{Q}}^0[\mathtt{W}'^c_b]=0$. Rewrite this null equation in term of $\mathbb{Q}$ and of $A$, we obtain $$
\mathbb{E}[\int_{[0,\infty]}\ind_{\mathtt{W}'^c_b}(v)dA_v]=0.
$$
Notice that $\mathtt{W}'_b\supset \mathtt{W}'_c$ for any pair $0\leq b\leq c<T$, and $\cup_{b<c<T}\mathtt{W}'_c=\mathtt{W}'_b$. The process $\ind_{\mathtt{W}'_b}\mathsf{p}_{b+}, b\in[0,T)$, has all its trajectories càdlàg and it is $\check{\mathbb{G}}^0$ adapted (i.e. $\check{\mathbb{G}}^0$ optional). 

Notice that, when we replace $\mathsf{p}_{+}$ by $\ind_{\mathtt{W}'}\mathsf{p}_{+}$, all the computations in the previous paragraphs remain valid. The theorem is proved with this function $\ind_{\mathtt{W}'}\mathsf{p}_{+}$. \ok

\bcor
Assume the same condition as in Theorem \ref{acp}. The process $(\mathsf{p}_{u+}(\omega,u): u\in[0,T))$ is $\mathbb{F}$ progressively measurable on $[0,T)$. The process $\int_{[0,t]}\mathsf{p}_{v+}(v)dA_v, t\in[0,T)$, is equal to the $\mathbb{F}$ optional dual projection of the increasing process $\ind_{[\tau,\infty)}$ on $[0,T)$.
\ecor

\textbf{Proof.} By the right continuity in $t\in\mathbb{R}_+$ of the process $\mathsf{p}_{t+}$, for any $(\omega,u)\in\Omega\times[0,\infty]$, we can write$$
\mathsf{p}_{u+}(\omega,u)=\lim_{N\uparrow \infty}\sum_{k=1}^{\lfloor 2^NT\rfloor}\mathsf{p}_{v^N_k+}(\omega,u)\ind_{\{v^N_{k-1}\leq u<v^N_k\}},\ t\in[0,T),
$$ 
where $v_k=\frac{k}{2^N}$. We notice that, for any $k\leq \lfloor 2^NT\rfloor$, for any $a\in\mathbb{R}$, for any $t\in[0,T)$, the set $$
\{(\omega,u)\in\Omega\times[0,t]:\ \mathsf{p}_{v^N_k+}(\omega,u)\ind_{\{v^N_{k-1}\leq u<v^N_k\}}\leq a\ \}
$$ 
is in $$
\dcb
\{v^N_{k-1}\leq \hat{\tau}<v^N_k\}\cap\{\hat{\tau}\leq t\}\cap\check{\mathcal{G}}_{v^N_{k}}
+\{v^N_{k-1}\leq \hat{\tau}<v^N_k\}^c\cap\{\hat{\tau}\leq t\}\cap\{\emptyset, \Omega\times[0,\infty]\}
\subset\mathcal{F}_{t+\frac{2}{2^N}}\otimes\mathcal{B}[0,t].
\dce
$$
Consequently, for any $a\in\mathbb{R}$, for any $t\in[0,T)$, the set $$
\{(\omega,u)\in\Omega\times[0,t]:\ \mathsf{p}_{u+}(\omega,u)\leq a\ \}
\in \mathcal{F}_{t+\epsilon}\otimes\mathcal{B}[0,t],
$$ 
for any $\epsilon>0$. Now applying \cite[Chapitre IV, n$^\circ$14]{DM}, we prove that the process $(\mathsf{p}_{u+}(\omega,u): u\in[0,T))$ is $\mathbb{F}$ progressively measurable on $[0,T)$.

Let $K$ be a bounded $\mathbb{F}$ optional process. Then the process $\check{K}_t(\omega,u)=K_t(\pi(\omega,u))$, for $t\in\mathbb{R}_+,(\omega,u)\in\Omega\times[0,\infty]$, is $\check{\mathbb{F}}$ optional (which can be proved by a usual monotone class theorem argument). Applying Theorem \ref{acp}, for $a\in[0,T)$,
$$
\dcb
\mathbb{E}[K_\tau\ind_{\{\tau\leq a\}}]
&=\mathbb{E}[K_{\tau\wedge a}\ind_{\{\tau\leq \tau\wedge a\}}]
=\check{\mathbb{E}}[\check{K}_{\hat{\tau}\wedge a}\ind_{\{\hat{\tau}\leq \hat{\tau}\wedge a\}}]\\
&=\check{\mathbb{E}}^0[\check{K}_{\hat{\tau}\wedge a}\ind_{\{\hat{\tau}\leq \hat{\tau}\wedge a\}}\mathsf{p}_{(\hat{\tau}\wedge a)+}]
=\check{\mathbb{E}}^0[\check{K}_{\hat{\tau}}\ind_{\{\hat{\tau}\leq a\}}\mathsf{p}_{\hat{\tau}+}]
=\mathbb{E}[\int_{[0,a]}K_v\mathsf{p}_{v+}(v)dA_v].
\dce
$$
The process $\int_{[0,t]}\mathsf{p}_{v+}(v)dA_v, t\in[0,T)$, is finite càdlàg and $\mathbb{F}$ adapted. This proves the theorem. \ok

\bcor\label{disint}
Assume the same condition as in Theorem \ref{acp}. For any $b\in[0,T)$, for any bounded $\check{\mathcal{G}}_b$ measurable function $h$, the process $(\int_{[0,b]}h(u)\mathsf{p}_{t+}(u)dA_u, t\in[b,T))$ is a càdlàg version of the $(\mathbb{Q},\mathbb{F})$ uniformly integrable martingale $\mathbb{E}[h(\tau)\ind_{\{\tau\leq b\}}|\mathcal{F}_t], t\in[b,T)$. 
\ecor

\brem
This lemma is much more precise than Lemma \ref{Mint}.
\erem

\textbf{Proof.}
For any $\mathbb{F}$ stopping time $S\in[b,T)$, according to Theorem \ref{acp}, we write$$
\dcb
\mathbb{E}[h(\phi)\ind_{\{\tau\leq b\}}]
&=\check{\mathbb{E}}[h\ind_{\{\hat{\tau}\leq b\}}]
=
\check{\mathbb{E}}^0[h\ind_{\{\hat{\tau}\leq b\}}\mathsf{p}_{\check{S}+}]
&=\mathbb{E}[\int_{[0,b]}h(u)\mathsf{p}_{S+}(u)dA_u].
\dce
$$
Notice that the process $\ind_{\{\hat{\tau}\leq b\}}h\ind_{[b,T)}\mathsf{p}_{+}$ is $\check{\mathbb{G}}$ optional. According to Lemma \ref{oopp}, the process $\int_{[0,b]}h(u)\mathsf{p}_{t+}(u)dA_u, t\in[b,T)$, is $\mathbb{F}$ optional. The above identities together with \cite[Theorem 4.40]{Yan} implies that this process is a càdlàg $(\mathbb{Q},\mathbb{F})$ uniformly integrable martingale on $[b,T)$. On the other hand, by the monotone class theorem, using Theorem \ref{acp} property 3, it can be checked that $$
\int_{[0,b]}h(u)\mathsf{p}_{t+}(u)dA_u
=\mathbb{E}[h(\tau)\ind_{\{\tau\leq b\}}|\mathcal{F}_t]. \ \ok
$$

\

\bl\label{nullset}
Assume the same condition as in Theorem \ref{acp}. Let $b\in[0,T)$. Then, $$
\int_{[0,b]} \ind_{\{\mathsf{p}_{b+}(u)=0\}}\mathsf{p}_{t+}(u) dA_u=0,
$$
for any $t\in(b,T)$.
\el

\textbf{Proof.} 
We note that, according to  Corollary \ref{disint}, $\int_{[0,b]} \ind_{\{\mathsf{p}_{b+}(u)=0\}}\mathsf{p}_{t+}(u) dA_u$ is a càdlàg non negative $(\mathbb{Q},\mathbb{F})$ martingale on $[0,T)$ and$$
\int_{[0,b]} \ind_{\{\mathsf{p}_{b+}(u)=0\}}\mathsf{p}_{t+}(u) dA_u
=\mathbb{E}[\ind_{\{\mathsf{p}_{b+}(\tau)=0\}}\ind_{\{\tau\leq b\}}|\mathcal{F}_t].
$$
Taking the expectation$$
\mathbb{E}[\ind_{\{\mathsf{p}_{b+}(\tau)=0\}}\ind_{\{\tau\leq b\}}]
=\check{\mathbb{E}}[\ind_{\{\mathsf{p}_{b+}=0\}}\ind_{\{\hat{\tau}\leq b\}}]
=\check{\mathbb{E}}^0[\ind_{\{\mathsf{p}_{b+}=0\}}\ind_{\{\hat{\tau}\leq b\}}\mathsf{p}_{b+}]=0.\ \ok
$$

There is a last corollary.

\bcor
Let $A$ be a non negative $\mathbb{F}$ adapted increasing càdlàg process such that $A_\infty=1$ and $A_{t}< 1$ for $t\in\mathbb{R}_+$. Let $(M^u: u\in[0,\infty])$ be an i$\! M(\mathbb{Q},\mathbb{F})$ family. Let $\check{\mathbb{Q}}^0$ be the Cox measure on the product space associated with $(\check{\mathbb{Q}}|_{\check{\mathcal{F}}_\infty},\check{\mathbb{F}},\check{A},\hat{\tau})$. Then, $(M^u: u\in[0,\infty])$ is differentiable on $[0,T)$ with respect to $A$, if and only if $\check{\mathbb{Q}}$ is absolutely continuous with respect to $\check{\mathbb{Q}}^0$ on $\check{\mathcal{G}}_t$ for all $t\in[0,T)$.

\ecor

\textbf{Proof.} The condition is necessary by Theorem \ref{acp}. Suppose that the "if" condition holds. For $t\in[0,T)$ let $\xi$ denote the density function $\frac{d\check{\mathbb{Q}}}{d\check{\mathbb{Q}}^0}$ on $\check{\mathcal{G}}_t$. There exists a $\mathcal{F}_t\otimes\mathcal{B}[0,\infty]$ measurable function $\mathsf{p}(\omega,u)$ such that $$
\mathsf{p}(\hat{\tau}\nmid t)=\check{\mathbb{E}}^0[\xi|\check{\mathcal{F}}_t\vee\sigma(\hat{\tau}\nmid t)].
$$ 
Then, for $B\in\mathcal{F}_t$, for $0\leq u\leq t$, $$
\dcb
\mathbb{E}[\ind_B M^u_t]
&=\mathbb{E}[\ind_B\ind_{\{\tau\leq u\}}]
=\check{\mathbb{E}}[\ind_B(\pi)\ind_{\{\hat{\tau}\leq u\}}]\\
&=\check{\mathbb{E}}^0[\ind_B(\pi)\ind_{\{\hat{\tau}\leq u\}}\xi]
=\check{\mathbb{E}}^0[\ind_B(\pi)\ind_{\{\hat{\tau}\leq u\}}\mathsf{p}(\hat{\tau}\nmid t)]
=\mathbb{E}[\ind_B\int_{[0,u]}\mathsf{p}(u)dA_u],
\dce
$$
where the last equality is obtained by the definition of $\check{\mathbb{Q}}^0$. \ok

\

\section{Model defined by a stochastic differential equation}\label{sdetime}

Theorem \ref{acp} says that a model with a differentiable i$\!M$ is absolutely continuous with respect to a Cox measure model. As Cox model can be easily computed, a model with differentiable i$\! M_Z$ can easily be handled by Girsanov's theorem. One may wonder if necessary to introduce the notion of differentiable i$\! M_Z$ models. There are at least two reasons to introduce the differentiable models. Firstly, a lot of model are defined directly with an i$\! M_Z$. The notion of the differentiability gives us an effective method to check if the model is absolutely continuous with respect to the Cox measure. Secondly, the theoretical approach by Girsanov's theorem may fail to be helpful for practice purpose (cf. \cite{w}) for the lack of the computability of the density function. Variable models are required to meet the needs in practice. 

We present in this section a class of models which are defined through their dynamic equation and which have differentiable i$\!M$. These results are established in \cite{songmodel}.

\subsection{i$\! M_Z$}\label{imz}

Let $Z$ to be a $(\mathbb{Q},\mathbb{F})$ supermartingale such that $0\leq Z\leq 1$. We call such a $Z$ an Azéma supermartingale. We notice that for any i$\! M$, the process $(M^u_u)_{u\in\mathbb{R}_+}$ is an Azéma supermartingale. On the other hand, in market modeling, data calibrated from the real market can be represented by an Azéma's supermartingale $Z$. An important question is, therefore, if there exists a model such that $(M^u_u)_{u\in\mathbb{R}_+}$ coincides with $Z$. We introduce the following definition (cf. \cite{JS2}).

\bd
An increasing family of positive martingales issued from $1-Z$ (in short i$\!M_Z(\mathbb{P},\mathbb{F})$ or simply i$\!M_Z$) is an i$\!M$ family $(M^u: u\in[0,\infty])$ which satisfies the following conditions: for any $0\leq u\leq t<\infty$, $M^u_u=1-Z_u$ and $M^u_t \leq  1-Z_t$.
\ed

The theorem below is an immediate consequence of the Theorem \ref{theoremIM}.

\bethe \label{martingaleM} 
Let $(M^u: u\in[0,\infty])$ to be an i$M(\mathbb{Q},\mathbb{F})$ family associated with the probability measure $\mathbb{Q}$. Let $Z$ be an Azéma's supermartingale. Then, $(M^u: u\in[0,\infty])$ is an i$\!M_Z$ family if and only if $\mathbb{Q}[t<\tau|\mathcal{F}_t]=Z_t$ for $t\geq 0$.
\ethe

\

\subsection{$\natural$-equation and $\natural$-pair}

We suppose the condition:

\textbf{Hy}(Z): $1-Z_t<0, 1-Z_{t-}<0$ for $t\in(0,\infty)$. 

Let $Z=M-A$ be the $(\mathbb{Q},\mathbb{F})$ canonical decomposition of $Z$ with $M$ a $(\mathbb{Q},\mathbb{F})$ local martingale and $A$ a non-decreasing $\mathbb{F}$ predictable process. Notice that ${^{\mathbb{F}\cdot p}}(1-Z)_t=1-Z_{t-}+\Delta_tA>0$ for any $0<t<\infty$, where the superscript ${^{\mathbb{F}\cdot p}}$ denotes the $(\mathbb{Q},\mathbb{F})$ predictable projection. 
We define, for $0<u<\infty$, $\widetilde{m}^u_t=\int_u^t \frac{-dM_s}{{^{\mathbb{F}\cdot p}}(1-Z)_s}$. Since obviously $d\widetilde{m}^u_t=d\widetilde{m}^v_t$ for $0<u<v\leq t<\infty$, we omit the superscripts and we denote simply
$$
d\widetilde{m}_t=\frac{-dM_t}{{^{\mathbb{F}\cdot p}}(1-Z)_t},\ t\in(0,\infty).
$$
Let $\mathbb{D}$ design the space of all càdlàg $\mathbb{F}$ adapted processes. Let $m>0$ be an integer. Let $\mathbf{Y}=(Y_1,\ldots,Y_m)$ be an $m$-dimensional $(\mathbb{P},\mathbb{F})$ local martingale, and $\mathbf{F}=(F_1,\ldots,F_m)$ be a Lipschitz functional from $\mathbb{D}$ into the set of $m$-dimensional locally bounded $\mathbb{F}$ predictable processes in the sense of \cite{protter}. For $0< u<\infty$, for any ${\mathcal {F}}_u$-measurable random variable $x$, we consider the stochastic differential equation determined by the pair $(\mathbf{F},\mathbf{Y})$: $$
(\natural_u)  \left\{\dcb
dX_t&=&X_{t-} d\widetilde{m}_t+\mathbf{F}(X)_{t}^\top d\mathbf{Y}_t,\ t\in[u,\infty),\\
X_u&=&x.
\dce
\right.
$$

We will call the pair $(\mathbf{F},\mathbf{Y})$ a $\natural$-pair if it satisfies the following conditions, for any $1\leq j\leq m$, for any $u>0$ and for any $X,X'\in\mathbb{D}$:
\ebe
\item[(i)]
The process $t\in[u,\infty)\rightarrow\frac{F_j(X)_t}{{^{\mathbb{F}\cdot p}}(1-Z)_t-X_{t-}}\ind_{\{(1-Z_{t-})-X_{t-}\neq 0\}}$ is integrable with respect to $Y_j$, and satisfies the inequality:
$$
\Delta_t\tilde{m}-\frac{1}{{^{\mathbb{F}\cdot p}}(1-Z)_t-X_{t-}}\ind_{\{{^{\mathbb{F}\cdot p}}(1-Z)_t-X_{t-}\neq 0\}}\mathbf{F}(X)_t^\top\Delta_t\mathbf{Y}>-1, \  t\in[u,\infty).
$$

\item[(ii)]
The process $t\in[u,\infty)\rightarrow\frac{F_j(X)_t}{X_{t-}}\ind_{\{X_{t-}\neq 0\}}$ is integrable with respect to $Y_j$, and satisfies the inequality:
$$
\Delta_t\tilde{m}+\frac{1}{X_{t-}}\ind_{\{X_{t-}\neq 0\}}\mathbf{F}(X)_t^\top\Delta_t\mathbf{Y}\geq -1, \  t\in[u,\infty).
$$

\item[(iii)]
If $0\leq X,X'\leq 1$, the process $t\in[u,\infty)\rightarrow\frac{F_j(X)_t-F_j(X')_t}{X_{t-}-X'_{t-}}\ind_{\{X_{t-}-X'_{t-}\neq 0\}}$ is integrable with respect to $Y_j$, and satisfies the inequality:
$$
\Delta_t\tilde{m}+\frac{1}{X_{t-}-X'_{t-}}\ind_{\{X_{t-}-X'_{t-}\neq 0\}}(\mathbf{F}(X)_t-\mathbf{F}(X')_t)^\top\Delta_t\mathbf{Y}\geq -1, \  t\in[u,\infty).
$$
\dbe

The following theorem proves that the set of $\natural$-pairs is not empty.

\bethe\label{markov}
Let $g(t,x)$ be any bounded continuously differentiable function defined on $\mathbb{R}_+\times \mathbb{R}$ taking values in $\mathbb{R}^m$. Let $\varphi$ be a $C^\infty$ increasing function on $\mathbb{R}_+$ such that $|\varphi(x)|\leq 2$ and $|\frac{\varphi(x)}{x}|\leq 1$. For $t\in \mathbb{R}_+$, we introduce the set $\mathsf{G}_{t}$ of $\mathbf{z}\in\mathbb{R}^m$ satisfying the two conditions:$$
\dcb
\circ&:&2\left|g(t,x)^\top\mathbf{z}\right|<1+\Delta_t\tilde{m},\mbox{ for $x\in\mathbb{R}$,}\\
\\
\circ\circ&:&\left[-\varphi'({^{\mathbb{F}\cdot p}}(1-Z_{t-})-x)\varphi(x)g(t,x)+\varphi({^{\mathbb{F}\cdot p}}(1-Z_{t-})-x)\varphi'(x)g(t,x)\right.\\
&&\hspace{3cm}+\left.\varphi({^{\mathbb{F}\cdot p}}(1-Z_{t-})-x)\varphi(x)g'(t,x)\right]^\top\mathbf{z}
>-(1+\Delta_t\tilde{m}),\mbox{ for $x\in\mathbb{R}$}.
\dce
$$
(Here $g'(t,x)$ denotes the derivative with respect to $x$.) Then, for any $t\in\mathbb{R}_+$, the random set $\mathtt{G}_t$ is not empty, and the set-valued process $\mathtt{G}$ is $\mathbb{F}$ optional. There exists an $m$-dimensional $\mathbb{F}$ local martingale $\mathbf{Y}=(Y_1,\ldots,Y_m)$ whose jump at $t\in\mathbb{R}_+$, if it exists, is contained in $\mathsf{G}_{t}$. Let $$
\mathbf{F}(X)_t=f(t,X_{t-})=\varphi({^{\mathbb{F}\cdot p}}(1-Z)_t-X_{t-})\varphi(X_{t-}) g(t,X_{t-}), \ X\in\mathbb{D}.
$$
Then, the above conditions (i), (ii) and (iii) with strict inequality $>-1$ instead of $\geq -1$ are satisfied for the pair $(\mathbf{F},\mathbf{Y})$.
\ethe

\

\subsection{Model defined by the $\natural$-equation}

The $\natural$-equation defines an i$\!M_Z$.

\bethe Let $(\mathbf{F},\mathbf{Y})$ be a $\natural$-pair. For $ 0<u<\infty$, consider the equation $(\natural_u)$ associated with $(\mathbf{F},\mathbf{Y})$. Let $(L^u_t:t\in[u,\infty))$ denote the solution of the
equation $(\natural_u)$ with the initial condition $L^u_u=1-Z_u$.
Set $L^u_\infty=\lim_{t\rightarrow\infty}L^u_t$. Set $$ \dcb
M^u_u=(1-Z_u),\\
M^u_t=\inf_{v\in\mathtt{Q},u<v\leq t}(L^v_t)^+\wedge(1-Z_t),\ t\in(u,\infty].\\
\dce
$$
Set finally$$
\dcb
M^0_t&=&\inf_{u\in\mathbb{Q}, 0<u\leq t}M^u_t,\ t\in(0,\infty],\\
M^0_0&=&\lim_{t\downarrow 0}M^0_t \ \mbox{ (which exists)},\\
M^\infty_t&=&1,\ \mbox{ for $t\in[0,\infty]$}.
\dce
$$
Then, for $0<u<\infty$, $M^u$ is $\mathbb{P}$ indistinguishable to $L^u$ on $[u,\infty]$, and $(M^u:0\leq u\leq \infty)$ is an i$M_Z$. 
\ethe 

The above i$\!M_Z$ will be said to be associated with the $\natural$-equation as well as the probability measure $\mathbb{Q}^\natural$ constructed in Theorem \ref{theoremIM} with this i$\!M_Z$ will be said to be associated with the $\natural$-equation.

\

\subsection{The differentiability}

Let $g$ be a $C^\infty$ function on $\mathbb{R}$ with a compact support. Consider the $\natural$-equation associated with the following $\natural$-pair $(\mathbf{F},\mathbf{Y})$ of the type in Theorem \ref{markov}:$$
\dcb
\Delta_tY&\in&\mathtt{G}_t,\\
\mathbf{F}(X)_t&=&\varphi({^{\mathbb{F}\cdot p}}(1-Z)_t-X_{t-})\varphi(X_{t-}) g(X_{t-}), \ X\in\mathbb{D}.
\dce
$$ 
We suppose moreover that $\varphi(x)=x$, for $x\in[0,1]$. Consider the i$\!M_Z=(M^u:0\leq u\leq \infty)$ associated with $(\mathbf{F},\mathbf{Y})$. Since $0\leq M^u\leq 1-Z$, we have$$
\varphi({^{\mathbb{F}\cdot p}}(1-Z)-M^u_-)\varphi(M^u_-)=({^{\mathbb{F}\cdot p}}(1-Z)-M^u_-)M^u_-.
$$
This means that $M^u$ satisfies the following stochastic differential equation $$
\left\{\dcb
dX_t=X_{t-}d\tilde{m}_t+({^{\mathbb{F}\cdot p}}(1-Z)_t-X_{t-})X_{t-}g(X_{t-})^\top d\mathbf{Y}_t,\ u\leq t<\infty,\\
X_u=x,
\dce
\right.
$$
with $M^u_u=1-Z_u$. Theorem 39 and Theorem 65 in \cite[Chapter V Section 10]{protter} are applicable to such an equation. Let $x\rightarrow \Xi_t^{u}(x)$ be the associated stochastic differential flow.

\bethe\label{regularity}
Let $0< t<\infty$. Let$$
\kappa_v=(1+\Delta_v\widetilde{m}-(1-Z_{v-})g(1-Z_{v-})^\top\Delta_v\mathbf{Y}), \ \mbox{ for $v\in(0,t]$}.
$$
Then, when $\Delta_vA>0$, $$
M^v_t-M^{v-}_t=\Xi_t^v(1-Z_v)-\Xi_t^v(1-Z_v-\kappa_v\Delta_vA),\ v\in(0,t].
$$
When $\Delta_vA=0$,$$
\dcb
\lim_{u\uparrow v}\frac{M^v_t-M^u_t}{A_v-A_u}&=\frac{d\Xi_t^v}{dx}(1-Z_v)\kappa_v,\ &\mbox{ for $v\in(0,t]$},\\
\\
\lim_{u\downarrow v}\frac{M^v_t-M^u_t}{A_v-A_u}&=\frac{d\Xi_t^v}{dx}(1-Z_v),\ &\mbox{ for $v\in(0,t)$}.
\dce
$$
Consequently, the i$\!M_Z$ is differentiable with respect to $A$ on $[0,\infty)$ with the density function$$
\mathsf{p}_t(u)=\frac{d\Xi_t^u}{dx}(1-Z_u)1\!\!1_{\{\Delta_uA=0\}}
+\frac{\Xi_t^u(1-Z_u)-\Xi_t^u(1-Z_u-\kappa_u\Delta_uA)}{\Delta_uA}1\!\!1_{\{\Delta_uA>0\}},\ u\in(0,t], t\in\mathbb{R}_+.
$$
\ethe

\

\section{Copulas and the ordering statistics of random times}\label{copula}

In defaultable market modeling, we need to take the order statistics of a family of default times. In this section we prove a sufficient condition which ensure that the order statistics satisfy the differentiability property for their i$\! M_Z$.

\subsection{Ordering of functions on $\{1,\ldots,k\}$}\label{ordering}

We begin with recalling the order statistics. Let $\mathfrak{a}$ be a function defined on $\{1,\ldots,k\}$  (where $k>0$ is an integer) taking values in $[0,\infty]$. Let $\{a_1,\ldots,a_k\}$ denote the values of $\mathfrak{a}$. Let $$
R^\mathfrak{a}(i)=R^{\{a_1,\ldots,a_k\}}(i)=\sum_{j=1}^k \ind_{\{a_j<a_i\}}+\sum_{j=1}^k \ind_{\{j<i,a_j=a_i\}}+1.
$$
The map $i\in\{1,\ldots,k\}\rightarrow R^\mathfrak{a}(i)\in\{1,\ldots,k\}$ is a bijection. Let $\rho^\mathfrak{a}$ be its inverse. Define $\uparrow\!\!\!\mathfrak{a}=\mathfrak{a}(\rho^\mathfrak{a})$. We check that $\uparrow\!\!\!\mathfrak{a}$ is a non decreasing function on $\{1,\ldots,k\}$ taking the same values of $\mathfrak{a}$. $\uparrow\!\!\!\mathfrak{a}(i)$ represents the $i$th smallest values among $\{a_1,\ldots, a_k\}$.

Let $k\in\mathbb{N}^*$ and $\tau_1,\ldots,\tau_k$ be $k$ random times. Consider the function $\mathfrak{t}$ on $\{1,\ldots,k\}$ taking respectively the values $\{\tau_1,\ldots,\tau_k\}$. We define $\sigma_{i}=\uparrow\!\!\!\mathfrak{t}(i), 1\leq i\leq k$. Note that, if the $\tau_i$ are stopping times with respect to some filtration, the $\sigma_{i}$ are stopping times with respect to the same filtration, because$$
\{\sigma_{i}\leq t\}=\cup_{I\subset\{1,\ldots,k\}, \sharp I=i}\{\tau_j\leq t, \forall j\in I\}, \ t\geq 0.
$$
This same equation shows that there exists a Borel function $\mathfrak{s}_{i}$ on $[0,\infty]^k$ such that $\sigma_{i}=\mathfrak{s}_{i}(\tau_1,\ldots,\tau_k)$.

\bl\label{inex}
Let $t\in[0,\infty]$. Let $1\leq i\leq k$. Denote $\mathfrak{S}_i=\{I\subset\{1,\ldots,k\}: \sharp I=i\}$. Also, for $I\subset \{1,\ldots,k\}$ denote $B_I=\{\tau_j\leq u, \forall j\in I\}$. We have $$
\mathbb{Q}[\sigma_i\leq u|\mathcal{F}_t]
=\sum_{S\subset \mathfrak{S}_i: S\neq \emptyset}(-1)^{1+|S|}\mathbb{Q}[\cap_{I\in S}B_I|\mathcal{F}_t].
$$
\el

\textbf{Proof.} We have$$
\dcb
\mathbb{Q}[\sigma_i\leq u|\mathcal{F}_t]
&=&\mathbb{Q}[\cup_{I\subset\{1,\ldots,k\}: \sharp I=i}\{\tau_j\leq u, \forall j\in I\}|\mathcal{F}_t]
=\mathbb{Q}[\cup_{I\in\mathfrak{S}_i}B_I|\mathcal{F}_t].
\dce
$$
By the inclusion-exclusion formula, the last term becomes$$
\dcb
\mathbb{Q}[\cup_{I\in\mathfrak{S}_i}B_I|\mathcal{F}_t]
&=&\sum_{S\subset \mathfrak{S}_i: S\neq \emptyset}(-1)^{1+|S|}\mathbb{Q}[\cap_{I\in S}B_I|\mathcal{F}_t]. \ \ok
\dce
$$

\subsection{Continuously differentiable copulas}

We refer to \cite{nelsen} for the notion of the copulas.

\bethe
Let $C(x_1,\ldots,x_k)$ be a continuously differentiable copulas. Let $\mathbf{M}^i=(M^{i,u}: u\in[0,\infty]), 1\leq i\leq k$, be i$\! M_Z$ satisfying the differentiability condition at a point $T\in[0,\infty]$ with respect to a same increasing process $A$. Let $\tau_i, 1\leq i\leq k$, be a family of random items such that their conditional law given $\mathcal{F}_T$ has the (multi-dimensional) distribution function $(u_1,\ldots,u_k)\rightarrow C(M^{1,u_1}_T,\ldots,M^{k,u_k}_T)$ on $0\leq u_1,\ldots,u_k\leq T$. Then, the order statistics of the $\tau_i, 1\leq i\leq k$, satisfy the differentiability condition for their i$\! M$ at every $t\in[0,T]$ with respect to $A$.
\ethe

\textbf{Proof.} For any $J\subset\{1,\ldots,k\}$, denote by $C_J$ the $J$-marginal copulas of $C$, which is also continuously differentiable. Let $\mathsf{p}^i_T$ be a density function at $T$ of $\mathbf{M}^i$. For $0\leq u\leq t\leq T$, we compute$$
\dcb
\mathbb{Q}[\tau_j\leq u, j\in J|\mathcal{F}_t]
&=&\mathbb{Q}[\mathbb{Q}[\tau_j\leq u, j\in J|\mathcal{F}_T]|\mathcal{F}_t]

=\mathbb{Q}[C_J(M^{1,u}_T,\ldots,M^{k,u}_T)|\mathcal{F}_t].
\dce
$$
We write the Ito's formula (cf. \cite[(2.52)]{jacod}) to the expression $C_J(M^{1,u}_T,\ldots,M^{k,u}_T)$ with respect to variable $u$.$$
\dcb
&&C_J(M^{1,u}_T,\ldots,M^{k,u}_T)\\
&=&C_J(M^{1,0}_T,\ldots,M^{k,0}_T)
+\sum_{j\in J}\int_0^u\frac{\partial C_j}{\partial x_j}(M^{1,u-}_T,\ldots,M^{k,u-}_T)\mathsf{p}^j_T(s)dA_s\\
&&+\sum_{s\in(0,u]}\ind_{\{\Delta_sA> 0\}}\left(\frac{C_J(M^{1,s}_T,\ldots,M^{k,s}_T)-C_J(M^{1,s-}_T,\ldots,M^{k,s-}_T)}{\Delta_s A}-\sum_{j\in J}\frac{\partial C_j}{\partial x_j}(M^{1,s-}_T,\ldots,M^{k,s-}_T)\mathsf{p}^j_T(s)\right)\Delta_sA\\
&=&\sum_{j\in J}\int_0^u\frac{\partial C_j}{\partial x_j}(M^{1,u-}_T,\ldots,M^{k,u-}_T)\mathsf{p}^j_T(s)dA^c_s\\
&&+\sum_{s\in(0,u]}\ind_{\{\Delta_sA> 0\}}\left(\frac{C_J(M^{1,s}_T,\ldots,M^{k,s}_T)-C_J(M^{1,s-}_T,\ldots,M^{k,s-}_T)}{\Delta_s A}\right)\Delta_sA,
\dce
$$
where $A^c$ denotes the continuous part of the process $A$. Hence, there exists a $\mathcal{F}_T\otimes\mathcal{B}[0,\infty]$ measurable function $\xi_J(\omega,u)$ such that$$
C_J(M^{1,u}_T,\ldots,M^{k,u}_T)
=\int_{[0,u]}\xi_J(s)dA^c_s
+\sum_{s\in[0,u]}\xi_J(s)\Delta_sA
=\int_{[0,u]}\xi_J(s)dA_s.
$$
Consequently,$$
\mathbb{Q}[\tau_j\leq u, j\in J|\mathcal{F}_t]
=\mathbb{Q}[\int_{[0,u]}\xi_J(s)dA_s|\mathcal{F}_t]
=\int_{[0,u]}\mathbb{Q}[\xi_J(s)|\mathcal{F}_t]dA_s.
$$
Now using the same notations in Lemma \ref{inex}, we write$$
\dcb
&&\mathbb{Q}[\sigma_i\leq u|\mathcal{F}_t]
=\mathbb{Q}[\cup_{I\in\mathfrak{S}_i}B_I|\mathcal{F}_t]\\
&=&\sum_{S\subset \mathfrak{S}_i: S\neq \emptyset}(-1)^{1+|S|}\mathbb{Q}[\cap_{I\in S}B_I|\mathcal{F}_t]\\

&=&\sum_{S\subset \mathfrak{S}_i: S\neq \emptyset}(-1)^{1+|S|}\mathbb{Q}[B_{\cup_{I\in S}I}|\mathcal{F}_t]\\

&=&\sum_{S\subset \mathfrak{S}_i: S\neq \emptyset}(-1)^{1+|S|}\int_{[0,u]}\mathbb{Q}[\xi_{\cup_{I\in S}I}(s)|\mathcal{F}_t]dA_s\\

&=&\int_{[0,u]}\left(\sum_{S\subset \mathfrak{S}_i: S\neq \emptyset}(-1)^{1+|S|}\mathbb{Q}[\xi_{\cup_{I\in S}I}(s)|\mathcal{F}_t]\right)dA_s.\ \ok
\dce
$$

\brem
The above theorem presents an invariant principle for the class of differentiable i$\! M$. We point out here another obvious invariant principle, i.e. the class of differentiable i$\! M$ is invariant by absolutely continuous change of probability measures.
\erem

\

\section{Conditional expectation}\label{cex}

In this section we fix a random time $\tau$ whose Azéma's supermartingale is $Z=M-A$ as in subsection \ref{imz}. We suppose that its i$\!M$ family is differentiable on $[0,T)$ with respect to $A$ ($A_\infty=1$) with a density function $\mathsf{p}$. This section is devoted to the computation of conditional expectations under this assumption.

\

\subsection{parametered optional projection}

Recall that in \cite[Proposition 3.]{SY}, for any non negative $\mathcal{B}(\mathbb{R}_+)\otimes\mathcal{A}\otimes\mathcal{B}[0,\infty]$ measurable function $H(t,\omega,u)$, there is constructed a non negative $\mathcal{O}(\mathbb{F})\otimes\mathcal{B}[0,\infty]$ measurable function, denoted by ${^o}H$ such that, for every $u\in[0,\infty]$, ${^o}H(u)$ is a version of the $(\mathbb{Q},\mathbb{F})$ optional projection of the process $H(u)$. In the following, we will introduce a variant of \cite{SY}'s notion. 

We introduce $\mathbf{C}$ the family of bounded $\mathcal{F}_\infty\otimes\mathcal{B}[0,\infty]$ measurable functions such that, for each member $F$ of that family, there exists an $\mathcal{O}(\mathbb{F})\otimes\mathcal{B}[0,\infty]$ measurable function, denoted by $({^o}\!F)_t(\omega,u)$, $(t,\omega,u)\in\mathbb{R}_+\times\Omega\times[0,\infty]$, which satisfies the following conditions:

\ebe
\item[$\flat$]
The process $({^o}\!F)_t(t), t\in\mathbb{R}_+$, is $\mathbb{F}$ optional. 
\item[$\flat\flat$]
For any $u\in[0,\infty]$, $({^o}\!F)_t(u), t\in\mathbb{R}_+$, is a version of the $(\mathbb{Q},\mathbb{F})$ optional projection of $F(u)$.
\item[$\flat\flat\flat$]
For any $\mathbb{F}$ stopping times $U$, the process $\ind_{\{U\leq t\}}({^o}\!F)_t(U), t\in\mathbb{R}_+$, is $\mathbb{F}$ optional. For any pair of $\mathbb{F}$ stopping times $U,S$ such that $U\leq S$, for any $B\in\mathcal{F}_S$,
$$
\mathbb{E}^0[\ind_B F(U)]
=\mathbb{E}[\ind_B\ ({^o}\!F)_S(U)].
$$
\dbe
Clearly, the family $\mathbf{C}$ is a linear space and it contains the function $F\equiv 1$. Notice that, if $F,F'$ are two members of the family $\mathbf{C}$ with $0\leq F\leq F'$, by section theorem (cf. \cite{Yan}), the process $\ind_{\{U\leq t\}}({^o}\!F)_t(U), t\in\mathbb{R}_+$, is bounded from below by zero and is overestimated by the process $\ind_{\{U\leq t\}}({^o}\!F')_t(U), t\in\mathbb{R}_+$. Let $F$ be a bounded non negative $\mathcal{F}_\infty\otimes\mathcal{B}[0,\infty]$ measurable function which is the increasing limit of a sequence $(F_n)_{n\in\mathbb{N}^*}$ in $\mathbf{C}$. For any $n\in\mathbb{N}$, let ${^o}\!F_n$ be a $\mathcal{O}(\mathbb{F})\otimes\mathcal{B}[0,\infty]$ measurable function satisfying the above conditions $\flat$ to $\flat\flat\flat$. with respect to $F_n$. We set ${^o}\!F=\sup_{n\in\mathbb{N}^*}{^o}\!F_n$. Then ${^o}\!F$ is $\mathcal{O}(\mathbb{F})\otimes\mathcal{B}[0,\infty]$ measurable function and it satisfies the conditions $\flat$ to $\flat\flat\flat$ with respect to $F$. This observation shows that the family $\mathbf{C}$ is a functional monotone class in the sense of \cite[Theorem 1.4]{Yan}. For any $C\in\mathcal{F}_\infty$ and $D\in\mathcal{B}[0,\infty]$, taking $({^o}\!F)_t(u)=\mathbb{E}[\ind_C|\mathcal{F}_t]\ind_D(u)$, we see that the family $\mathbf{C}$ contains equally all functions of the form $F=\ind_C\ind_D$. By the monotone class theorem (cf. \cite[Theorem 1.4]{Yan}), the family $\mathbf{C}$ contains all bounded $\mathcal{F}_\infty\otimes\mathcal{B}[0,\infty]$ measurable function. By a usual limit procedure we extend the definition of ${^o}\!F$ also to all non negative functions $F$.

We can now introduce the following definition.

\bd
For any bounded or non negative $\mathcal{F}_\infty\otimes\mathcal{B}[0,\infty]$ measurable function $F$, ${^o}\!F$ denote any $\mathcal{O}(\mathbb{F})\otimes\mathcal{B}[0,\infty]$ measurable function satisfying the above conditions $\flat$ to $\flat\flat\flat$.
\ed

\brem
Notice that, if the measurable space $(\Omega,\mathcal{A})$ has nice topological property, the function ${^o}\!F$ can be defined more explicitly with the prediction process in \cite{aldous}.
\erem

The next lemma is a direct consequence of the definition of ${^o}\!F$.

\bl\label{poF}
Let $F$ be any bounded or non negative $\mathcal{F}_\infty\otimes\mathcal{B}[0,\infty]$ measurable function and $U$ be a $\mathbb{F}$ stopping time. If $F(U)$ is integrable, the maps ${^o}F_t(U)$, $t\in[U,\infty)$, define a $(\mathbb{Q},\mathbb{F})$ uniformly integrable càdlàg martingale on $[U,\infty)$. Consequently, the $\mathbb{F}$ predictable projection of $F(U)$ on $[U,\infty)$ is given by ${^o}F_{t-}(U), t\in[U,\infty)$.
\el

\bl\label{e0}
For any non negative non decreasing $\mathbb{F}$ adapted càdlàg process $A$, for any bounded or non negative $\mathcal{F}_\infty\otimes\mathcal{B}[0,\infty]$ function $F$, for any version of ${^o}\!F$, we have$$
\mathbb{E}[\int_{[0,t]} f(u)\ F(u)dA_u]
=\mathbb{E}[\int_{[0,t]} f(u)\ {^o}\!F_t(u)dA_u],
$$
for any $t\in\mathbb{R}_+$ and any bounded $\cap_{s>t}(\mathcal{F}_s\otimes\mathcal{B}[0,\infty])$ measurable function $f$. 
\el

\textbf{Proof.} We give a proof when $F\geq 0$. Let $c(s)=\inf\{t\in\mathbb{R}_+: A_t>s\}$. Using the results in subsection \ref{changevariable}, we write
$$
\dcb
&&\mathbb{E}[\int_{[0,t]} f(u)\ F(u)dA_u]\\

&=&\mathbb{E}[\int_{[0,\infty)}\ind_{\{c(s)\leq t\}} f(c(s))\ F(c(s)) ds]\\

&=&\int_{[0,\infty)}\mathbb{E}[\ind_{\{c(s)\leq t\}} f(c(s))\ F(c(s))]\ ds, \ \mbox{ by Fubini's theorem,}\\

&=&\int_{[0,\infty)}\mathbb{E}[\ind_{\{c(s)\leq t\}} f(c(s))\ ({^o}\!F)_t(c(s))]\ ds \\ &&\mbox{ because $c(s)$ is a $\mathbb{F}$ stopping time and $\ind_{\{c(s)\leq t\}}f(c(s))$ is $\mathcal{F}_t$ measurable,}\\

&=&\mathbb{E}[\int_{[0,t]} f(u)\ {^o}\!F_t(u)dA_u].\  \ok
\dce
$$

\bl\label{mart0}
In the same setting as in the previous lemma, we have$$
\check{\mathbb{E}}^0[F|\check{\mathcal{G}}^0_t]\ind_{\{\hat{\tau}\leq t\}}
={^o}F_t\ind_{\{\hat{\tau}\leq t\}}.
$$
\el

\textbf{Proof.} We rewrite the computation in the above lemma in a different form:
$$
\dcb
\check{\mathbb{E}}^0[f\ind_{\{\hat{\tau}\leq t\}}F]
&=&\mathbb{E}[\int_{[0,t]} f(u)\ F(u)dA_u]\\

&=&\mathbb{E}[\int_{[0,t]} f(u)\ {^o}\!F_t(u)dA_u]

=\mathbb{E}^0[f\ind_{\{\hat{\tau}\leq t\}}\ {^o}\!F_t]. \ \ok
\dce
$$

\

\subsection{Formulas of conditional expectation}

In this subsection we will show that the conditional expectation with respect to $\mathcal{G}_t$ can be computed in term of the density function $\mathsf{p}$ and in term of the conditional expectation with respect to $\mathcal{F}_t$.

The following lemma has been called "key lemma" in the literature (cf. \cite[Lemma 3.1.2]{BR}). Recall $Z_t=\mathbb{E}[\ind_{\{t<\tau\}}|\mathcal{F}_t]$.

\bl\label{keylemma}
For any bounded $\mathtt{H}\in\mathcal{G}_\infty$, 
$$
\mathbb{E}[\mathtt{H}|\mathcal{G}_t]\ind_{\{t<\tau\}}
=\frac{\mathbb{E}[\mathtt{H}\ind_{\{t<\tau\}}|\mathcal{F}_t]}{Z_t}\ind_{\{t<\tau\}},\ t\in[0,\infty).
$$
\el

\textbf{Proof.}
Firstly, let $h$ be a bounded $\mathcal{F}_t\otimes\mathcal{B}[0,\infty]$ measurable function. We have$$
\dcb
\mathbb{E}[h(\tau\nmid t)\ind_{\{t<\tau\}}\mathtt{H}]

&=\mathbb{E}[h(\infty)\ind_{\{t<\tau\}}\mathtt{H}]

=\mathbb{E}[h(\infty)\mathbb{E}[\mathtt{H}\ind_{\{t<\tau\}}|\mathcal{F}_t]]\\
\\
&=\mathbb{E}[h(\infty)\ind_{\{t<\tau\}}\frac{\mathbb{E}[\mathtt{H}\ind_{\{t<\tau\}}|\mathcal{F}_t]}{\mathbb{E}[\ind_{\{t<\tau\}}|\mathcal{F}_t]}]

=\mathbb{E}[h(\tau\nmid t)\ind_{\{t<\tau\}}\frac{\mathbb{E}[\mathtt{H}\ind_{\{t<\tau\}}|\mathcal{F}_t]}{Z_t}\ind_{\{t<\tau\}}].
\dce
$$
We note that every quantities in this formula are bounded. Hence, we can take the right limit in $t\in\mathbb{R}_+$ and conclude that the same relationship remains valid when $h(\tau\nmid t)$ is replaced by bounded $\mathcal{G}_t$ measurable functions. \ok

\bethe\label{conditioning}
Let $b\in[0,T)$. For any non negative $\check{\mathcal{G}}^0_b$ measurable function $f$, for $t\in[0,b]$,
$$
\mathbb{E}[f(\phi)|\mathcal{G}_t]
=\frac{\mathbb{E}[f(\phi)\ind_{\{t<\tau\}}|\mathcal{F}_t]}{Z_t}\ind_{\{t<\tau\}}
+\frac{{^o}\!(f\mathsf{p}_{b+})_t(\tau)}{\mathsf{p}_{t+}(\tau)}\ind_{\{\mathsf{p}_{t+}(\tau)>0\}}\ind_{\{\tau\leq t\}}.
$$
\ethe

\textbf{Proof.} The formula on $\{t<\tau\}$ is the consequence of Lemma \ref{keylemma}. Consider the formula on $\{\tau\leq t\}$. Notice$$
\dcb
\mathbb{E}[f(\phi)|\mathcal{G}_t]\ind_{\{\tau\leq t\}}
=\left(\check{\mathbb{E}}[f|\check{\mathcal{G}}_t]\ind_{\{\hat{\tau}\leq t\}}\right)(\phi).
\dce
$$
Apply Theorem \ref{acp}. $$
\dcb
\check{\mathbb{E}}[f|\check{\mathcal{G}}_t]\ind_{\{\hat{\tau}\leq t\}}
&=&\frac{1}
{\mathsf{P}_t}\check{\mathbb{E}}^0[f\ind_{\{\hat{\tau}\leq b\}}\mathsf{P}_b|\check{\mathcal{G}}_t]\ind_{\{\mathsf{P}_{t}>0\}}\ind_{\{\hat{\tau}\leq t\}}

=\frac{1}
{\mathsf{p}_{t+}}\check{\mathbb{E}}^0[f\mathsf{p}_{b+}|\check{\mathcal{G}}_t]\ind_{\{\mathsf{p}_{t+}>0\}}\ind_{\{\hat{\tau}\leq t\}}\\

&=&\frac{{^o}\!(f\mathsf{p}_{b+})_t}{\mathsf{p}_{t+}}\ind_{\{\mathsf{p}_{t+}>0\}}\ind_{\{\tau\leq t\}},
\dce
$$
according to Lemma \ref{mart0}. \ok

\bcor\label{pts}
Let $b\in[0,T)$ and $t\in[0,b]$. We have $
\ind_{\{\tau\leq t\}}\mathsf{p}_{t+}(\tau)=\ind_{\{\tau\leq t\}}{^o}\!(\mathsf{p}_{b+})_t(\tau),
$
and $\ind_{\{\mathsf{p}_{t+}(\tau)>0\}}\ind_{\{\tau\leq t\}}=\ind_{\{\tau\leq t\}}$, $\mathbb{Q}$ almost surely. Moreover, ${^o}\!(\mathsf{p}_{b+})_t(\omega,u)$ is another version of the density function at $t$.
\ecor

\textbf{Proof.} If we take $f\equiv 1$ in the formula of Theorem \ref{conditioning}, we get$$
\ind_{\{\tau\leq t\}}
=\frac{{^o}\!(\ind_{[0,b]}\mathsf{p}_{b+})_t(\tau)}{\mathsf{p}_{t+}(\tau)}\ind_{\{\mathsf{p}_{t+}(\tau)>0\}}\ind_{\{\tau\leq t\}}.
$$
This proves the first part of the lemma. The second part is the consequence of the identity$$
\mathbb{E}[h(\tau)\ind_{\{\tau\leq t\}}]
=\mathbb{E}[\int_{[0,t]} h(u)\ {^o}(\mathsf{p}_{b+})_t(u)dA_u].
$$
for any bounded $\check{\mathcal{G}}^0_t$ measurable function $h$, consequence of Lemma \ref{e0}. \ok

\brem
Notice that the new version ${^o}\!(\mathsf{p}_{b+})_t(u)$ of the density function has a simpler measurability than $\mathsf{p}_{t+}$. In fact, ${^o}\!(\mathsf{p}_{b+})_t(u)$ is $\mathcal{F}_t\otimes\mathcal{B}[0,\infty]$ measurable, for $t\in[0,b]$. An additional remark is that, for every $u\in[0,\infty]$, $t\rightarrow {^o}\!(\mathsf{p}_{b+})_t(u)$ is a $(\mathbb{Q},\mathbb{F})$ uniformly integrable martingale on $t\in[0,b]$. 
\erem

\

\section{The optional splitting formula}

We assume the same assumption as in Section \ref{cex}. In this section we consider the notion of the optional splitting formula introduced in \cite{songsplitting}. We have the following result.

\bethe\label{optionalsplitting}
The optional splitting formula on $[0,T)$ holds in the filtration $\mathbb{G}$ under the probability measure $\mathbb{Q}$. This means that, for any $\mathbb{G}$ optional process $X$, there exist a $\mathbb{F}$ optional process $X'$ and a $\mathcal{O}(\mathbb{F})\otimes\mathcal{B}[0,\infty]$ measurable function $X''$ such that the following $\mathbb{Q}$ indistinguishable identity holds:$$
X\ind_{[0,T)}=X'\ind_{[0,\tau)}\ind_{[0,T)}+X''(\tau)\ind_{[\tau,\infty)}\ind_{[0,T)}.
$$
\ethe

\textbf{Proof.} Let $X$ be a $\mathbb{G}$ optional process. By \cite{songsplitting} we know that there exists always a $\mathbb{F}$ optional process $X'$ such that $X\ind_{[0,\tau)}=X'\ind_{[0,\tau)}$. We need only to prove the optional splitting formula on the interval $[\tau,\infty)\cap[0,T)$.

According to \cite[Theorem 4.36]{Yan} there exists a $\mathbb{G}^0$ optional process $X^0$ which is $\mathbb{Q}$ indistinguishable from $X$. According to Lemma \ref{optionalinclusion} and \cite[Theorem 1.5]{Yan}, there exists a $\check{\mathbb{G}}$ optional process $\check{X}$ such that $X^0=\check{X}(\phi)$.

Let $\check{\mathbb{Q}}^0$ be the Cox measure on the product space associated with the triplet $(\check{\mathbb{Q}}|_{\check{\mathcal{A}}},\hat{\tau},\check{A})$. Under the Cox measure $\check{\mathbb{Q}}^0$, the hypothesis$(H)$ (cf. \cite[Lemma 4.2.1]{BJR}) is satisfied so that the optional splitting formula holds (cf. \cite{songsplitting}). Hence, the following $\check{\mathbb{Q}}^0$ indistinguishable identity holds:$$
\check{X}=\check{X}'\ind_{[0,\hat{\tau})}+\check{X}''(\hat{\tau})\ind_{[\hat{\tau},\infty)},
$$
for a $\check{\mathbb{F}}$ optional process $\check{X}'$ and a $\mathcal{O}(\check{\mathbb{F}})\otimes\mathcal{B}[0,\infty]$ measurable function $\check{X}''$. 

Notice that, according to Theorem \ref{acp}, $\check{\mathbb{Q}}$ is absolutely continuous with respect to $\check{\mathbb{Q}}^0$ on the set $\check{\mathcal{G}}_b$ for any $b\in[0,T)$. Hence, 
$$
\check{X}\ind_{[\hat{\tau},T)}=\check{X}''(\hat{\tau})\ind_{[\hat{\tau},T)}
$$
also is a $\check{\mathbb{Q}}$ indistinguishable identity. The map $\check{X}''$ is a function of four variables : $$
\check{X}''_t((\omega,u),x),\ \mbox{with } (t,(\omega,u))\in\mathbb{R}_+\times(\Omega\times[0,\infty]), x\in [0,\infty].
$$ 
We define $X''_t(\omega,x)=\check{X}''_t(\phi(\omega),x)$ which is a $\mathcal{O}(\mathbb{F})\otimes\mathcal{B}[0,\infty]$ measurable function. Then, we have the following $\mathbb{Q}$ indistinguishable identity:
$$
X^0\ind_{[\tau,T)}=X''(\tau)\ind_{[\tau,T)}. \ok
$$

\bcor
For any $t\in[0,T)$, $\mathcal{F}_t\vee\sigma(\tau\nmid t)$ completed by $(\mathbb{Q},\mathcal{G}_\infty)$ null sets is equal to $\mathcal{G}_t$. Moreover, $\sigma(H_\tau:\mbox{ $H$ a $\mathbb{F}$ optional process})$ completed by $(\mathbb{Q},\mathcal{G}_\infty)$ null sets is equal to $\mathcal{G}_\tau$.
\ecor

\textbf{Proof.} It is the consequence of the preceding theorem and \cite[Theorem 3.4 and 3.6]{songsplitting}. \ok 

\brem
The property proved in the above corollary is important in regard to \cite{BR} where such a property under the name condition \textbf{G.1} and \textbf{G.2} is required to establish results. 
\erem

\

\section{Enlargement of filtration formula}

It is essential, for a market model based on the progressive enlargement of filtration to be useful, to know the hypothesis$(H')$ holds, i.e. to know if all $\mathbb{F}$ martingale $X$ remains $\mathbb{G}$ semimartingale (otherwise there will exist arbitrage). If it is the case, it is important to know the semimartingale decomposition of $X$ in $\mathbb{G}$. This section is devoted to that question. 

We assume the same assumption as in Section \ref{cex}. We recall that $Z$ is the Azéma supermartingale of $\tau$ and $M$ is its martingale part (in $\mathbb{F}$). We begin with a well-known result in \cite{JY, JY2, J}.

\bl\label{Jlemma} 
Let $X$ be a bounded $(\mathbb{Q},\mathbb{F})$ martingale. Let $B^X$ the $(\mathbb{Q},\mathbb{F})$ predictable dual projection of the jump
process $t\rightarrow \Delta X_\tau\ind_{\{0<\tau\leq t\}}$. Then, $$
X_{\cdot\wedge \tau} - \int_0^{\cdot\wedge
\tau}\frac{1}{Z_{s-}}(d\cro{M,X}_s+dB^X_s)
$$ is a $(\mathbb{Q},\mathbb{G})$  local martingale.
\el

\textbf{Notation.} Let $\mathsf{P}$ be the $(\check{\mathbb{Q}}^0,\check{\mathbb{G}})$ martingale on the time interval $[0,T)$ introduced in Theorem \ref{acp}. Let $X$ be a
bounded $(\mathbb{Q},\mathbb{F})$ martingale. Then the bounded process $\check{X}$ (where $\check{X}=X(\pi)$) is a $(\check{\mathbb{Q}}^0,\check{\mathbb{F}})$ martingale. By \cite{BJR} the hypothesis$(H)$ is satisfied between $\check{\mathbb{F}}$ and $\check{\mathbb{G}}$ under the Cox measure $\check{\mathbb{Q}}^0$. It results that $\check{X}$ is a $(\check{\mathbb{Q}}^0,\check{\mathbb{G}})$ martingale. We can therefore compute their $(\check{\mathbb{Q}}^0,\check{\mathbb{G}})$ predictable bracket process on the time interval $[0,T)$, denoted by $\cro{\check{X},\mathsf{P}}$. Using the predictable splitting formula (see \cite[Lemme(4.4)]{J}), there exist a $\mathcal{P}(\check{\mathbb{F}})$ measurable process $Y'$ and a $\mathcal{P}(\check{\mathbb{F}})\otimes\mathcal{B}[0,\infty]$ measurable function $Y''(u)$ such that$$
\cro{\check{X},\mathsf{P}}\ind_{[0,T)}
=Y'\ind_{[0,\hat{\tau}]}\ind_{[0,T)}+Y''(\hat{\tau})\ind_{(\hat{\tau},\infty)}\ind_{[0,T)}.
$$
Recall that $\ind_{[0,T)}\mathsf{p}_{+}=\ind_{[\hat{\tau},T)}\mathsf{P}$. With an argument by the monotone class theorem we see that there exists a $\mathcal{P}(\mathbb{F})\otimes\mathcal{B}[0,\infty]$ measurable function $\cro{X,\mathsf{p}_{+}}_t(\omega,s), (t,\omega)\in\mathbb{R}_+\times\Omega, s\in[0,\infty]$, such that $Y''_t((\omega,u),s)=\cro{X,\mathsf{p}_{+}}_t(\pi(\omega,u),s)$.

\bethe
Under the assumption of this section, for any bounded $(\mathbb{Q},\mathbb{F})$ martingale $X$, the process
$$
\dcb
X_{t}-X_0-\int_{0}^{t}\ind_{\{s\leq \tau\}}\frac{1}{Z_{s-}}(d\cro{M,X}_s+dB^X_s)

-\int_{0}^{t} \ind_{\{\tau<s\}}\frac{1}{\mathsf{p}_{s-}(\tau)}d\cro{X,\mathsf{p}_+}_s(\tau)
\dce
$$ is a $(\mathbb{Q},\mathbb{G})$ local martingale on the time interval $t\in [0,T)$, where $\mathsf{p}_{-}$ denotes the left limit process of $\mathsf{p}_{+}$.
\ethe

\textbf{Proof.} According to Lemma \ref{Jlemma}, $X^{\tau}$ is a $(\mathbb{Q},\mathbb{G})$ special semimartingale with the drift process  $$
\int_0^{t\wedge
\tau}\frac{1}{Z_{s-}}(d\cro{M,X}_s+dB^X_s),\ t\in\mathbb{R}_+.
$$ 
We have already indicated that $\check{X}$ is a bounded $(\check{\mathbb{Q}}^0,\check{\mathbb{F}})$ martingale. It is hence a $(\check{\mathbb{Q}}^0,\check{\mathbb{G}})$ martingale. By Girsanov's theorem, 
$$
\check{X}-\check{X}^{\hat{\tau}}
-\frac{1}{\mathsf{P}_{-}}\ind_{(0,T)}\centerdot\cro{\check{X}-\check{X}^{\hat{\tau}},\mathsf{P}}
$$
is a $(\check{\mathbb{Q}},\check{\mathbb{G}})$ local martingale on $[0,T)$. Note that, using the notation fixed previously to the theorem,$$
\dcb
\frac{1}{\mathsf{P}_{-}}\ind_{(0,T)}\centerdot\cro{\check{X}-\check{X}^{\hat{\tau}},\mathsf{P}}
=\frac{1}{\mathsf{P}_{-}}\ind_{(\hat{\tau},T)}\centerdot \cro{\check{X},\mathsf{P}}
=\frac{1}{\mathsf{p}_{-}}\ind_{(\hat{\tau},T)}\centerdot Y''(\hat{\tau}).
\dce
$$
Now we pull the above martingale property back to $\Omega$ by $\phi$ to conclude that 
$$
X-X^\tau
-\frac{1}{\mathsf{p}_{-}(\tau)}\ind_{(\tau,T)}\centerdot \cro{X,\mathsf{p}_+}(\tau)
$$
is a $(\mathbb{Q},\mathbb{G})$ local martingale on $[0,T)$. This proves the theorem because $X=X^\tau+(X-X^\tau)$. \ok

\

\textbf{Acknowledgment} This research benefited from the support of the "Chair Markets in Transition", under the aegis of Louis Bachelier laboratory, a joint initiative of Ecole polytechnique, Université d'Evry Val d'Essonne and Fédération Bancaire Française.

\

\end{document}